\documentclass{amsart}
\usepackage{amssymb}
\usepackage{amsfonts,dsfont}
\usepackage{color}
\usepackage{amsmath}
\usepackage{comment}
\usepackage{graphicx}
\usepackage{todonotes,bm}
\usepackage[bookmarks,bookmarksopen=false,                  pdfauthor={Autor},                  pdftitle={Titel},                  colorlinks=true,                  linkcolor=blue,                  citecolor=blue,                  urlcolor=blue,                ]{hyperref}
\allowdisplaybreaks
\usepackage{enumitem}

\newtheorem{theorem}{Theorem}

\newtheorem{lemma}[theorem]{Lemma}

\newtheorem{proposition}[theorem]{Proposition}
\newtheorem{remark}[theorem]{Remark}

\renewcommand{\proof}{\noindent\textbf{Proof:\ }}
\newcommand{\pbox}{\hfill$\Box$\\}
\newcommand{\V}{\mathbf{V}}
\newcommand{\R}{\mathbb{R}}

\DeclareMathOperator{\sinc}{sinc}

\definecolor{darkviolet}{rgb}{0.58,0,0.83} 

\definecolor{lblue}{rgb}{0.4,0.4,0.75}

\def\bd{\mathbf}
\newcommand{\dgt}{DGT }

\begin{document}

\title[Time-frequency analysis on flat tori]{Time-frequency analysis on flat tori and Gabor frames in finite dimensions}
\author{L. D. Abreu}
\address{Faculty of Mathematics \\
	University of Vienna \\
	Oskar-Morgenstern-Platz 1 \\
1090 Vienna, Austria}
\email{abreuluisdaniel@gmail.com}
\email{michael.speckbacher@univie.ac.at}
\author{P. Balazs}
\address{Acoustics Research Institute\\Austrian Academy of Sciences\\Wohllebengasse 12-14, 1040 Vienna\\ Austria }
\email{peter.balazs@oeaw.ac.at}
\email{nicki.holighaus@oeaw.ac.at}
\author{N. Holighaus}
\author{F. Luef }
\address{Department of Mathematics, NTNU Trondheim, 7041 Trondheim, Norway}
\email{franz.luef@ntnu.no}
\author{M. Speckbacher}

\subjclass[2010]{42C40, 46E15, 42C30, 46E22, 42C15}
\keywords{short-time Fourier transform, flat torus, \dgt frames, Feichtinger algebra, sampling theory}
\thanks{The authors would like to thank Hans Georg Feichtinger for valuable discussions and comments, to Antti Haimi for his input during the early stages of this work and to Karlheinz Gr{\"o}chenig for bringing \cite{bannert2013discretized} to our attention. This research was supported by the Austrian Science Fund (FWF) through the projects  P-31225-N32 (L.D.A.), P 34624 (P.B.) Y-1199, J-4254 (M.S.), as well as I 3067-N30 (N.H.). (L.D.A.) was also supported by Portuguese funds through CIDMA-Center for Research and Development in Mathematics and Applications, and FCT– \textquotedblleft Funda\c{c}\~{a}o para a Ci\^{e}ncia e
a Tecnologia\textquotedblright, within project UIDB/04106/2020 and
UIDP/04106/2020.}
\date{}
\maketitle

\begin{abstract}
\noindent We provide the foundations of a Hilbert space theory for the short-time
Fourier transform (STFT) where the flat tori
\begin{equation*}
\mathbb{T}_{N}^2=\mathbb{R}^2/(\mathbb{Z}\times N\mathbb{Z})=[0,1]\times \lbrack
0,N] 
\end{equation*}
act as phase spaces. We work on an $N$-dimensional subspace $S_{N}$ of distributions periodic in time and frequency in the dual $S_0'(\R)$ of the Feichtinger algebra $S_0(\R)$ and equip it with an inner product. To construct the Hilbert space $S_{N}$ 
we apply a suitable double periodization operator to $S_0(\mathbb{R})$. On $S_{N}$, the STFT is applied as the usual STFT defined on $S_0'(\R)$.
This STFT is a continuous extension of the finite discrete Gabor transform from the lattice onto the entire flat torus. As such, 
sampling theorems on flat tori lead
to Gabor frames in finite dimensions. For Gaussian windows, one is lead to spaces of analytic functions and the construction allows to prove a necessary and sufficient Nyquist rate type result, which is the analogue, for Gabor frames in finite dimensions, of a well known result of Lyubarskii and Seip-Wallst{\'e}n for Gabor frames with Gaussian windows and which, for $N$ odd, produces an explicit \emph{full spark Gabor frame}. 
The compactness of the phase space, the finite dimension of the signal spaces and our sampling theorem 
offer practical advantages in some applications. We illustrate this by discussing 
a problem of current research interest: recovering signals from the zeros of their noisy spectrograms. 
\end{abstract}

\section{Introduction}
The \emph{short-time Fourier transform (STFT)} is the central instrument of time-frequency analysis. The most classical setting considers the analysis of functions $f$ with respect to windows $g$, both contained in $\mathrm{L}^2(\mathbb{R})$, defined as
\begin{equation}
\mathbf{V}_{g}f(x,\xi )=\int_{\mathbb{R}}f(t)\overline{g(t-x)}e^{-2\pi i\xi t}dt = \left\langle f,\mathbf{M}_{\xi }\mathbf{T}_{x}g\right\rangle =\langle f, {\pi}(x,\xi)g\rangle, \label{Gabor}
\end{equation}%
where $\mathbf{T}_{x}f(t)=f(t-x)$, $\mathbf{M}_{\xi }f(t)=e^{2\pi i\xi t}f(t)$, and $\pi(x,\xi)=\mathbf{M}_\xi \mathbf{T}_x$ define the \emph{translation}, \emph{modulation} and \emph{time-frequency shift} operators, respectively. By interpreting the brackets
as a duality pairing, this definition also holds for pairs of test function and
distribution spaces, like the Schwartz space and tempered distributions  $\mathcal{S(\mathbb{R})},\mathcal{S}^{\prime }(%
\mathbb{R})$~\cite{Charly} and, in particular, the Feichtinger algebra $ {S}_{0}(\mathbb{R})$ and its dual ${S}_{0}^{\prime }(\mathbb{R})$~\cite{fe06}. 

In this paper, we  consider the STFT acting on the $N$-dimensional 
space $S_{N}$ of time and frequency periodic distributions in ${S}^{\prime }_0(\mathbb{R})$, see definitions in Section 2 and \cite[Chapter~16.3]{Gazeau} or \cite[Chapter~6]{coro}.
This will lead to new phase spaces for the joint time and frequency values: the \emph{flat tori} $\mathbb{T}_{N}^2=[0,1]\times \lbrack
0,N]$, providing, as we will show, \emph{a continuous extension of the coefficient space of the discrete Gabor transform (DGT)}. The space $S_N$
is isometrically isomorphic to $\mathbb{C}^N$ equipped with the Euclidean norm and can similarly be obtained by sampling and periodization of $S_0(\mathbb{R})$~\cite{Kaib05}. This connection implies that results for the STFT on $S_N$ have implications for the \dgt on $\mathbb{C}^N$ and vice-versa. However, we will demonstrate that the STFT on the distribution space $S_N$, embedded into $ {S}^{\prime }_0(\mathbb{R})$, has much stronger structural properties, similar to those enjoyed by the STFT on $\mathrm{L}^2(\mathbb{R})$. As a continuous phase space extension of the DGT, the STFT on flat tori provides a natural way of defining \emph{off-the-grid} values, offering flexibility in applications and the chance of using continuous variable methods in finite Gabor analysis. 

In the case of Gaussian windows, we obtain spaces of analytic functions. The resulting possibility of using analytic complex variable tools will allow us to prove a necessary and sufficient Nyquist rate type result for Gabor frames with Gaussian windows in finite dimensions, which can be seen as the finite-dimensional analogue of the celebrated result of Seip-Wallst{\'e}n \cite{sewa92} and Lyubarskii \cite{ly92} for Gabor frames with Gaussian windows. The sufficient condition provides theoretical support to numerical procedures for increasing grid resolution, due to the principle of stable reconstruction using frames above the Nyquist rate. As a step in the proof of the sufficient Nyquist rate, we show that the STFT of any signal in $S_{N}$ with Gaussian window has exactly $N$ zeros, thereby making precise and proving the claim in~\cite{gardner2006sparse}. 
This property of the toric STFT   comes in handy for the problem of detecting noisy spectrogram zeros \cite{bardenet2020zeros,escudero2022efficient,courbot2023sparse,flandrin2015}. The potential advantages of the toric setting for such applications will be further discussed in the final paragraph of the paper.

Our methods are innovative in the sense that they allow to obtain results for finite sequences merely as a byproduct of the theory on $S_N$. But it must be noted that the relation between the continuous STFT and the discrete Gabor transform has been studied by several authors over the last 30 years, in particular by Janssen~\cite{Jans97}, and later by Kaiblinger~\cite{Kaib05} and S\o{}ndergaard~\cite{sondergaard2007finite,so05}.
Where the works of Janssen and S\o{}ndergaard are concerned with the construction of discrete Gabor frames and dual windows from Gabor systems on $S_0(\mathbb{R})$, Kaiblinger's work is concerned with finite dimensional approximation of dual windows for Gabor frames on $S_0(\mathbb{R})$. 
The sampling-periodization duality of the Fourier transform, succintly expressed in (generalizations of) Poisson's summation formula and considered in many works, including~\cite{Ausl96,Cool67,jalu19,Kaib05}, is central to these contributions.
In essence, the transition between $S_0(\mathbb{R})$ and $\mathbb{C}^N$ is achieved by studying a composition of periodization operators $\mathbf P_{(1)},\mathbf P_{(2)}$ on certain intervals with forward and inverse Fourier transforms $\mathcal F_{(1)},\mathcal F_{(2)}^{-1}$ as $\mathcal F_{(2)}^{-1}\mathbf P_{(2)}\mathcal F_{(1)}\mathbf P_{(1)}$.
Here, $\mathcal F_{(1)}$ is the Fourier transform of $\mathrm{L}^2(\mathbb{R})$ of a finite interval and $\mathcal F_{(2)}^{-1}$ is the inverse discrete Fourier transform. From this angle, the central deviation of the present paper from these prior works is that we consider $\mathcal F_{(1)}$ and $\mathcal F_{(2)}^{-1}$ to be distributional Fourier transforms on ${S}_{0}^{\prime }(\mathbb{R})$, such that  $\mathcal F_{(2)}^{-1}\mathbf P_{(2)}\mathcal F_{(1)}\mathbf P_{(1)} f$,  $f\in S_0(\mathbb{R})$, yields a doubly periodic distribution in ${S}_{0}^{\prime }(\mathbb{R})$ instead of a finite sequence in $\mathbb{C}^N$, enabling the subsequent application of the STFT on ${S}_{0}^{\prime }(\mathbb{R})$ instead of the \dgt.

\section{Overview}

We consider functions, distributions, and finite sequences, denoted by lower case latin letters $f,g$, greek letters  $\phi,\psi$, and sans font latin letters $\mathsf{f},\mathsf{g}$, respectively. For the latter, the discrete nature of the domain of $\mathsf{f},\mathsf{g}$ is emphasized by using square brackets for indexing, e.g. $\mathsf{f}[l]$. Operators are denoted by upper case letters $\mathbf{V},\mathbf{\Sigma}$. Exceptions from this convention are time-frequency shifts $\pi$, the Jacobi theta function $\vartheta$, and the Fourier transform $\mathcal F$, for which we adopt established notation.

With the Gaussian window $h_0(t)= e^{-\pi t^2}$, the Feichtinger algebra $ {S}_{0}(\mathbb{R})$ \cite{FeichingerSegal,jakobsen} is the space
\begin{equation*}
\begin{split}
  {S}_{0}(\mathbb{R}) & := \big\{ f\in \mathrm{L}^2(\mathbb{R})\colon \mathbf{V}_{h_0}f \in \mathrm{L}^1(\mathbb{R}^2)\big\},\ \text{equipped with the norm}\\
  \|f\|_{{S}_0} & :=\int_{\R^2}|\mathbf{V}_{h_0}f(x,\xi)|dxd\xi = \| \mathbf{V}_{h_0} f\|_{\mathrm{L}^1(\mathbb{R}^2)}.
\end{split}
\end{equation*}
We define the space $S_{N}$ as the span of $\{\epsilon_n\}_{n=0}^{N-1}$, the sequence of periodic delta trains \cite{Gazeau} 
\begin{equation}\label{def:deltatrain}
\epsilon _{n}:=\sum_{k\in \mathbb{Z}}\delta _{\frac{n}{N}+k}\subset S_0^\prime(\mathbb{R}),\qquad
n=0,...,N-1,
\end{equation}
and will show that $S_N$ can be characterized as the image of $
S_{0}(\mathbb{R})$ under \emph{the double periodization operator}
\begin{equation}\label{eq:quasiper1}
\mathbf{\Sigma}_{N}f:=\sum_{k_{1},k_{2}\in \mathbb{Z}}\mathbf{M}_{Nk
_{2}}\mathbf{T}_{k_{1}}f=\sum_{k_{1},k_{2}\in \mathbb{Z}}e^{2\pi
iNk_{2}\, \bullet}\cdot f(\, \bullet -k_{1})\text{.}
\end{equation}
It can be directly observed that $\mathbf{V}_{g}(\mathbf{\Sigma} _{N}f)$ is quasiperiodic, i.e.
\begin{equation}\label{eq:quasiper}
 \begin{array}{l}
\mathbf{V}_{g}(\mathbf{\Sigma} _{N}f)(x+1,\xi
)  =e^{-2\pi i\xi}\mathbf{V}_{g}(\mathbf{\Sigma} _{N}f)(x,\xi ),
\\ 
 \mathbf{V}_{g}(\mathbf{\Sigma} _{N}f)(x,\xi+N
)=\mathbf{V}_{g}(\mathbf{\Sigma }_{N}f)(x,\xi ).
\end{array} 
\end{equation}
Thus,\ the phase spaces of $\mathbf{V}_{g}\circ
\mathbf{\Sigma} _{N}$\ are the \emph{flat tori} $\mathbb{T}_{N}^2=[0,1]\times \lbrack
0,N]$.
As we will see,  $\mathbf{V}_g\circ \mathbf{\Sigma}_N: {S}_0(\R)\rightarrow \mathrm{L}^{2}(\mathbb{T}_N^2)$ and $\mathbf{V}_g:S_N\rightarrow \mathrm{L}^{2}(\mathbb{T}_N^2)$ have the same range in phase-space. It will often be convenient to jump from one to the other representation to simplify proofs.

The STFT on $S_N$ naturally introduces the compact phase space $\mathbb{T}_{N}^2$ for time-frequency analysis on finite, $N$-dimensional Hilbert spaces. Thereby, it provides a continuous model that, by construction, eliminates the truncation, or alternatively aliasing, errors usually associated with the transition from the STFT on $\mathrm{L}^2(\mathbb{R})$ to numerical implementations by means of the finite Gabor transform. That is not to say that these errors are removed: They are instead separated from the continuous model to the double periodization operator $\mathbf{\Sigma}_N$, i.e., the mapping from ${S}_{0}(\mathbb{R})$ onto $S_N\subset {S}_{0}^{\prime}(\mathbb{R})$.

As discussed in~\cite{Jans97,Kaib05,so05} in a slightly different formal framework, the
composition $\mathbf{V}_{g}\circ\mathbf{ \Sigma} _{N}$\ relates to the \emph{discrete Gabor transform} (DGT) on $\mathbb{C}^N$,
defined as
\begin{equation*}
\mathbf{V}_{\mathsf{g}}\mathsf{f}\left[ k,l\right] =\sum\limits_{m=0}^{N-1}\mathsf{f}
[m]\overline{\mathsf{g}[m-k]}e^{\frac{-2\pi ilm}{N}},\quad \mathsf{f},\mathsf{g}\in \mathbb{C}^{N},
\end{equation*}
We will show that $\mathbf{V}_{g}$ maps $ S_N$ into $\mathrm{L}^2(\mathbb{T}_N^2)$ 
and that $\mathbf{V}_{g}$ can be viewed as a continuous extension of $\mathbf{V}_{\mathsf{g}}\colon \mathbb{C}^N \rightarrow \mathbb{C}^{N\times N}$ to $\mathbb{T}^2_{N}$ in the sense of the following result.

\begin{theorem}\label{thm:finite-gabor} Let $f,g\in S_{0}(\mathbb{R})$ and let $\mathsf{f_N} = \mathbf{P}_N f,\mathsf{g_N} = \mathbf{P}_N g$, with the periodization operator
\[
  \mathbf{P}_N \colon S_{0}(\mathbb{R}) \rightarrow \mathbb{C}^N,\quad \text{defined by}\quad f\mapsto \Big(\sum_{j\in\mathbb{Z}} f(n/N-j)\Big)_{n=0}^{N-1}.
\]
Then, for $k,l\in 0,....,N-1$,
\begin{equation*}
{\mathbf{V}}_{g}( \mathbf{\Sigma} _{N}f)\left( \frac{k}{N},l\right) =N^{-1}\cdot\V_{\mathsf{g_N}}\mathsf{f_N}[k,l].
\end{equation*}
\end{theorem}

\begin{figure}
    \centering
    \includegraphics[width=12.5cm,trim={1.5cm 0 2cm 0},clip]{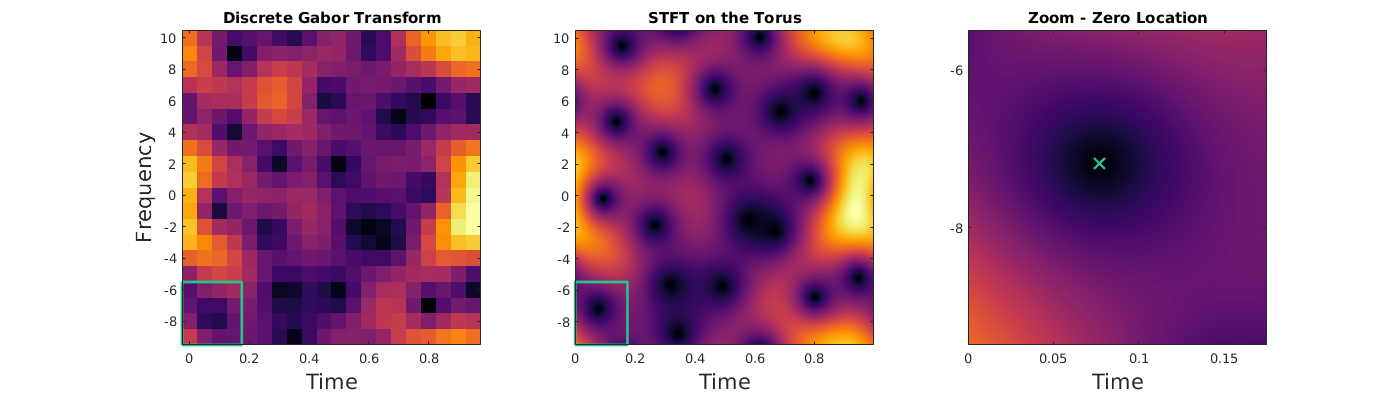}
    \caption{Comparison between the spectrograms obtained from the \dgt and the STFT on the torus, and their application for detecting zeros of the STFT: 
    Here, $\bd P_N f = \mathsf{f_N}$ was sampled from a complex Gaussian noise process $\mathcal{N}(0,I_N) + i\mathcal{N}(0,I_N)$, $N=20$. The right panel shows a zoom into the lower-left part of the spectrogram (indicated by the teal border), where the location of a spectrogram zero is marked.}
    \label{fig:torus_spectrograms}
\end{figure}

\noindent Specifically, the \dgt can be obtained by sampling the phase
space $\mathbb{T}^2_{N}$ of the continuous STFT restricted to $S_{N}$ on the grid points $(k/N,l)$,
providing a direct link between continuous and finite discrete
time-frequency analysis. An example spectrogram of both the \dgt and the STFT on the torus\footnote{A MATLAB script reproducing Figure \ref{fig:torus_spectrograms} can be downloaded from \url{https://www.oeaw.ac.at/en/ari/torusSTFT}.} is shown in Figure \ref{fig:torus_spectrograms}.

\begin{remark}\label{rem2}
  We chose the periods $(1,N)$ in Equation \eqref{eq:quasiper1} for notational convenience. An arbitrary choice of $(c,d)\in\mathbb{R}^2_+$ and $(a,b)\in\mathbb{R}^2$, with $cd = N$, leads to equivalent results on the phase space $(a,b)\,+\,[0,c]\times [0,d]$. When studying the approximation of the STFT by the \dgt, as in \cite{Kaib05}, it is usually more convenient to consider the symmetric convention  
   $(c,d)=(\sqrt{N},\sqrt{N})$ and $(a,b)=(-\sqrt{N}/2,-\sqrt{N}/2)$, such that an increase in $N$ symmetrically expands the considered phase space $\tilde{\mathbb{T}}_N^2 = [-\sqrt{N}/2,\sqrt{N}/2]\times [-\sqrt{N}/2,\sqrt{N}/2]$ and the sampling density within, covering all of $\mathbb{R}^2$ in the limit. The associated double  periodization operator is independent of $(a,b)$ and given by $\tilde{\mathbf{\Sigma}}_N\colon f \mapsto \sum_{k_1,k_2} e^{2\pi i \sqrt{N} k_2 \, \bullet}\cdot f(\, \bullet \, - \sqrt{N}k_1)$.

\end{remark}

\noindent We will study the Hilbert space
properties of the map 
\begin{equation*}
\mathbf{V}_{g} :S_{N}\rightarrow \mathrm{L}^{2}(\mathbb{T}_{N}^2)\text{,}
\end{equation*} 
and derive the Moyal-type orthogonality relation
$$
\int_{\mathbb{T}_N^2}\mathbf{V}_{g_1}\varphi_1(x,\xi)\overline{\mathbf{V}_ {g_2}\varphi_2(x,\xi)}dxd\xi=N\langle \varphi_1,\varphi_2\rangle_{S_N}\langle g_2,g_1\rangle_{\mathrm{L}^{2}},
$$
as well as inversion and
reproducing formulas similar to those of the continuous STFT. 

For the STFT
with dilated Gaussian windows $h_0^\lambda(t)=e^{-\pi \lambda t^{2}}$ we will obtain  a sampling theorem on the torus which leads to a full description of the frame set for finite  Gabor frames with Gaussian windows in $\mathbb{C}^{N}$. The proof uses a Bargmann-type transform, (which up to a weight is the STFT with $h_0^\lambda$ 
) whose action on the space $S_{N}$ has previously been considered in a slightly different form in \cite{LV}. Finally, combining the sampling theorem on the torus with Theorem~\ref{thm:finite-gabor}, we are lead to the following full description of the frame set for finite Gabor expansions using a periodized, dilated Gaussian window. As far as we could check, this is a completely new result.

\begin{theorem}
\label{cor:char-frame-set}
Let $\lambda>0$, $\mathsf{h_N^\lambda}=\mathbf{P}_Nh_0^\lambda$, and $\{(j_{k},l_{k})\}_{k=1,...,K}$ be a collection of distinct pairs of integers  $j_{k},l_{k}\in 0,....,N-1$. The following are equivalent:
\begin{enumerate}
    \item[1).] 
The set $
\{(j_{k},l_{k})\}_{k=1,...,K}$    gives
rise to a finite Gabor frame with window $\mathsf{h_N^\lambda}$, i.e., there are constants $A,B>0$ such that, for
every $\mathbf{f}\in \mathbb{C}^N$,   
\begin{equation*}
A\left\|  \mathsf{f} \right\|_{\mathbb{C}^{N}}^{2}\leq
\sum_{k=1}^{K}\left\vert \V_{\mathsf{h_{N}^\lambda}}\mathsf{f}[j_k,l_k]\right\vert ^{2}\leq B \left\|\mathsf{f}\right\| _{\mathbb{C}^{N}}^{2}\text{.}
\end{equation*}
\item[2).] One of the three following conditions is satisfied:
\begin{enumerate}
    \item[(i)] $N^{2}\geq K\geq N+1,$
    \item[(ii)] $K=N$ is odd, 
\item[(iii)] $K=N$ is even and $\sum_{k=1}^N(j_k,l_k)\notin N\mathbb{N}\times N\mathbb{N}. $
\end{enumerate}
\end{enumerate}
\end{theorem}

We emphasize that this result has been possible to prove only thanks to our Hilbert space theory for the STFT on flat tori, and that it strongly depends on the use of complex variable methods for quasi periodic analytic functions. This reinforces the suggestion that 
time-frequency analysis on the torus provides a rich theory which encompasses the theory of finite Gabor frames and leads to new insights, potential in applications and proof of results which were out of reach without the toric phase space. 

The paper is organized as follows. Some required properties of the Hilbert
space $S_{N}$ and the operator $\mathbf{\Sigma} _{N}$ are presented in Section \ref{sec:doublyperiodic}. Section~\ref{sec:main} contains the proof of
Theorem~\ref{thm:finite-gabor} above, explicit computations with dilated Gaussian windows  $h_0^\lambda$, and derivations of the Moyal-type formula, together with the inversion and reproducing kernel formulas. In the
last section, the window is specialized to be the Gaussian. The resulting Bargmann-type transform is defined, and several properties of its range space of entire functions with periodic constraints (a toric analogue of the Fock space) are studied in detail. All these properties are then used in the proof of the main result of the section: the sampling theorem on the torus. Finally, combining this result with Theorem~\ref{thm:finite-gabor}, we derive a full characterization of    finite Gabor frames with periodized and sampled Gaussian windows.

\section{Properties of \texorpdfstring{$S_{N}$}{SN} and  \texorpdfstring{$\mathbf{\Sigma }_{N}$}{SigmaN}}\label{sec:doublyperiodic}

\subsection{The Hilbert space \texorpdfstring{$S_{N}$}{SN} of time-frequency periodic
distributions}

The space $S_N$ appears in theoretical physics in coherent state approaches \cite{coro,Gazeau}. 
 By definition of $S_N$ it is clear that the family $\{\epsilon_n\}_{n=0}^{N-1}$ defined in \eqref{def:deltatrain}
forms a basis. Therefore, expanding $\varphi ,\psi \in S_{N}$ with respect to this basis
\begin{equation*}
\varphi =\sum_{n=0}^{N-1}a_{n}\epsilon _{n},\qquad \psi
=\sum_{n=0}^{N-1}b_{n}\epsilon _{n},
\end{equation*}
we can define an inner product on $S_N$ by 
\begin{equation*}
\left\langle \varphi ,\psi \right\rangle _{S_{N}}=\sum_{n=0}^{N-1}a_{n}%
\overline{b_{n}}.
\end{equation*}
Clearly, $S_{N}$ can be identified with $\mathbb{C}^{N}$ equipped with
the standard inner product, and $\{\epsilon _{n}\}_{n=0}^{N-1}$ forms an orthonormal basis of $S_{N}$ as
\begin{equation*}
\left\langle \epsilon _{n},\epsilon _{m}\right\rangle _{S_{N}}=\delta _{n,m}
\text{.}
\end{equation*}
Note that for every $\varphi\in S_N$
\begin{equation}\label{eq:double-per-S_N}
\mathbf{T}_{1}\varphi
=\varphi,\quad\text{ and }\quad \mathbf{M}_{N}\varphi =\varphi.
\end{equation}
Therefore, $S_N$ is a space of
  distributions that are {periodic in time and  frequency}. Actually,   $S_N$ contains all distributions in $S_0^\prime(\R)$ that satisfy \eqref{eq:double-per-S_N}, see \cite[page 262,
(16.12)]{Gazeau}.

\subsection{The double periodization operator \texorpdfstring{$\mathbf{\Sigma}_N$}{SigmaN}}

\noindent We formally define the double periodization operator as
\begin{equation*}
f\mapsto \mathbf{\Sigma} _{N}f=\sum_{k_{1},k_{2}\in \mathbb{Z}
}\mathbf{M}_{Nk_{2}}\mathbf{T}_{k_{1}}f.
\end{equation*}
The next lemma shows when and in which sense this object is well-defined.

\begin{lemma}
\label{lem:double-period-2} The operator $\mathbf{\Sigma} _{N}$ is well-defined from ${S}_0(\R)$ into ${S}_0'(\R)$ with unconditional weak-$\ast$ convergence in   ${S}_0^\prime(\mathbb{R})$.
\end{lemma}

\proof
If $f,g\in {S}_{0}(\mathbb{R})$, then 
\begin{align*}
|\mathbf{V}_{g}(\mathbf{\Sigma }_{N}f)(x,\xi )| &\leq \sum_{k_{1},k_{2}\in 
\mathbb{Z}}|e^{-2\pi ik_{1}\xi }\mathbf{V}_{g}f(x-k_{1},\xi -Nk
_{2})|
 \\
&\leq \sum_{k_{1},k_{2}\in \mathbb{Z}}|\mathbf{V}_{g}(\mathbf{T}_{-x}\mathbf{M}_{-\xi
}f)(k_{1},Nk_{2})| 
\\
&\leq C\Vert \mathbf{T}_{-x}\mathbf{M}_{-\xi }f\Vert _{S_{0}}\Vert g\Vert _{S_{0}}
 \\
&\leq C\Vert f\Vert _{S_{0}}\Vert g\Vert _{S_{0}}\leq \infty ,
\end{align*}
by \cite[Lemma 3.1.3 and Corollary 12.1.12]{Charly}. This implies that $\mathbf{\Sigma}
_{N}f\in {S}_{0}^{\prime }(\mathbb{R})$ is well-defined. If we choose $x=\xi=0$, then absolute weak-$\ast$ convergence of the series $ \mathbf{\Sigma} _{N}f$ in $S_0'(\R)$ follows which in turn implies unconditional convergence.
\pbox

{
\begin{remark}
It is easy to see that \[0 < \|\mathbf{V}_{g}(\mathbf{\Sigma}_{N}f) - \mathbf{V}_{g}f\vert_{\mathbb{T}^2_N}\|_{\mathrm{L}^1(\mathbb{T}^2_N)} \leq \|(1-\mathds{1}_{\mathbb{T}^2_N})\mathbf{V}_{g}f\|_{\mathrm{L}^1(\mathbb{R}^2)},\]for all $0\neq f\in S_0(\mathbb{R})$. That is, the double periodization introduces a nonzero, but controllable distortion in the STFT. Using the alternative phase space and double periodization convention introduced in Remark~\ref{rem2}, we obtain the analogous statement $0 < \|\mathbf{V}_{g}(\tilde{\mathbf{\Sigma}}_{N}f) - \mathbf{V}_{g}f\vert_{\tilde{\mathbb{T}}^2_N}\|_{\mathrm{L}^1(\tilde{\mathbb{T}}^2_N)} \leq \|(1-\mathds{1}_{\tilde{\mathbb{T}}^2_N})\mathbf{V}_{g}f\|_{\mathrm{L}^1(\mathbb{R}^2)}$, where 
\[
  \|(1-\mathds{1}_{\tilde{\mathbb{T}}^2_N})\mathbf{V}_{g}f\|_{\mathrm{L}^1(\mathbb{R}^2)} \overset{N\rightarrow \infty}{\longrightarrow} 0, 
\]
i.e., the distortion can be made arbitrarily small by increasing $N$. If $f$ is even in a weighted $\mathrm{L}^1$-type modulation space, i.e., $\omega\cdot \mathbf{V}_{g} f\in \mathrm{L}^1(\mathbb{R}^2)$ for some weight function $\omega \colon \mathbb{R}^2 \rightarrow [1,\infty)$, then we can quantify the rate of convergence. Using the equality $\mathbf{V}_{g}f(x,\xi) = \mathbf{V}_{\bd D g} \bd D f (\sqrt{N}x,\xi/\sqrt{N})$, where  $\bd D = \bd D_{\sqrt{N}}$ denotes the unitary dilation $\bd D_{\sqrt{N}} \colon f \mapsto N^{1/4}f(\sqrt{N}\, \bullet\, )$, we can derive a similar, but slightly less intuitive, statement on $\mathbb{T}^2_N$.
%
%
%
\end{remark}
}

\noindent As $\mathbf{\Sigma} _{N}f$ is periodic in time and frequency, we  can expand it
with respect to the orthonormal basis $\left\{ \epsilon _{n}\right\} _{n=0}^{N-1}$.
In the next lemma, this expansion is obtained
explicitly. We will also show that $\mathbf{\Sigma}_N$ is surjective as a mapping from $S_0(\mathbb{R})$ to $S_N$. To do so, we need to define the periodization operator $\mathbf{P}f(t)=\sum_{k\in\mathbb{Z}}f(t-k)$.

\begin{lemma}
\label{lem:double-period-1} For every $f\in {S}_{0}(\mathbb{R})$
\begin{equation}
\mathbf{\Sigma} _{N}f=\frac{1}{N} \mathbf{P}f\cdot \sum_{k\in \mathbb{Z}}\delta _{\frac{k}{N}} 
=\frac{1}{N}\sum_{n=0}^{N-1} \mathbf{P}f\left( \frac{n}{N}\right)
\epsilon _{n}\in S_N\text{,}  \label{eq:a_n}
\end{equation}
and 
\begin{equation}
\langle \mathbf{\Sigma} _{N}f,g\rangle _{{S}_{0}^{\prime } \times 
{S}_{0} }=\frac{1}{N}\sum_{n=0}^{N-1} \mathbf{P}f\left( 
\frac{n}{N}\right) \overline{ \mathbf{P}g\left( \frac{n}{N}\right) }=N\cdot
\langle \mathbf{\Sigma} _{N}f,\mathbf{\Sigma }_{N}g\rangle _{S_{N}}\text{.}
\label{eq:sigma_N-f,g}
\end{equation}
Moreover, $\mathbf{\Sigma}_N:{S}_0(\mathbb{R})\rightarrow S_N$ is surjective.
\end{lemma}

\proof
Let $f,g\in {S}_{0}(\mathbb{R})$. Since the Poisson summation formula holds for functions in ${S}_0(\mathbb{R})$ (see e.g. \cite[Corollary~12.1.5]{Charly}), we have that the following equality holds in the distributional sense
\begin{equation*}
\sum_{k\in \mathbb{Z}}e^{2\pi iktN}=\frac{1}{N}\sum_{k\in \mathbb{Z}}\delta
_{\frac{k}{N}}(t)\in S_{N}. \label{eq:exponential-sum}
\end{equation*}
This shows that  
\begin{equation*}
\mathbf{\Sigma} _{N}f =\sum_{k\in \mathbb{Z}}\mathbf{M}_{Nk}\sum_{l\in \mathbb{Z}}\mathbf{T}_{l}f=\sum_{k\in \mathbb{
Z}}\mathbf{M}_{Nk} \mathbf{P}f =\frac{1}{N}\sum_{k\in \mathbb{Z}}\delta _{
\frac{k}{N}}  \mathbf{P}f 
\end{equation*}
with unconditional weak-$\ast$ convergence   in ${S}_0^\prime(\mathbb{R})$. Hence,  
\begin{equation*}
\langle \mathbf{\Sigma }_{N}f,g\rangle _{{S}_{0}^{\prime } \times 
{S}_{0} }=\frac{1}{N}\left\langle \sum_{k\in \mathbb{Z}
}\delta _{\frac{k}{N}} \mathbf{P}f,g\right\rangle _{{S}_{0}^{\prime } \times {S}_{0} }.
\end{equation*}
 Let us write 
\begin{equation*}
\sum_{k\in \mathbb{Z}}\delta _{\frac{k}{N}}=\sum_{n=0}^{N-1}\sum_{k\in 
\mathbb{Z}}\delta _{\frac{n}{N}+k}=\sum_{n=0}^{N-1}\epsilon _{n}
\end{equation*}
where the change of summation order is justified by e.g. \cite[Corollary 12.1.5]{Charly}. Using the periodicity of $\mathbf{P}f$ then
yields 
\begin{align*}
\langle \mathbf{\Sigma} _{N}f,g\rangle _{{S}_{0}^{\prime }\times 
{S}_{0}} &=\frac{1}{N}\sum_{n=0}^{N-1}\left\langle
\sum_{k\in \mathbb{Z}}\delta _{\frac{n}{N}+k} \mathbf{P}f,g\right\rangle _{
{S}_{0}^{\prime }\times {S}_{0}} 
\\
&=\frac{1}{N}\sum_{n=0}^{N-1}\sum_{k\in \mathbb{Z}} \mathbf{P}f\left( \frac{n
}{N}+k\right) \overline{g\left( \frac{n}{N}+k\right) } 
\\
&=\frac{1}{N}\sum_{n=0}^{N-1} \mathbf{P}f\left( \frac{n}{N}\right)
\sum_{k\in \mathbb{Z}}\overline{g\left( \frac{n}{N}+k\right) } 
\\
&=\frac{1}{N}\sum_{n=0}^{N-1} \mathbf{P}f\left( \frac{n}{N}\right) \langle
\epsilon _{n},g\rangle _{{S}_{0}^{\prime }\times 
{S}_{0}}.
\end{align*} 
Hence, (\ref{eq:a_n}) holds. The first equality of (\ref{eq:sigma_N-f,g})
follows from the second to last equality above. Finally, the second equality of (\ref{eq:sigma_N-f,g}) results
from combining (\ref{eq:a_n}) and the first equality of (\ref{eq:sigma_N-f,g}).

 It thus remains to show
that $\mathbf{\Sigma} _{n}:{S}_{0}(\mathbb{R})\rightarrow S_{N}$ is
surjective. By \eqref{eq:a_n} it suffices to show that
there exists a family of function $f_n\in {S}_{0}(\mathbb{R})$, $n=0,\ldots,N-1$, satisfying $%
\mathbf{P} f_n\left( \frac{k}{N}\right) =\delta _{n}(k),$   $k=0,\ldots ,N-1$. Such  
functions obviously exists. Take for instance $f_n(t):=\sinc(Nt-n)\cdot e^{-\pi (t-n/N)^{2}}$ which is even a Schwartz function. \pbox

\section{Time-frequency analysis on flat tori}\label{sec:main}

\subsection{Basic properties of \texorpdfstring{$\mathbf{V}_g$}{Vgf} on \texorpdfstring{$S_N$}{SN}}

The STFT defined on $S_N$ is, as we we subsequently show, closely connected to the  \emph{Zak transform} which is defined as
\begin{equation*}
\mathbf{Z}f(x,\xi )=\sum_{k\in \mathbb{Z}}f(x- k)e^{2\pi i k\xi
}.
\end{equation*}
For later reference we state here some elementary facts about the Zak transform (see e.g. \cite{Charly}): \medskip
 
\noindent \emph{Quasiperiodicity}: 
\begin{equation}
\mathbf{Z}f(x,\xi +k)=\mathbf{Z}f(x,\xi ),\quad \text{and}
\quad \mathbf{Z}f(x+k ,\xi )=e^{2\pi i k\xi }\mathbf{Z}(x,\xi ),\quad k\in\mathbb{Z},
\label{eq:quasiperiodic}
\end{equation} \vspace{-0.3cm}

\noindent \emph{Action on time-frequency shifts: }
\begin{equation}
\mathbf{Z}(\mathbf{M}_{\omega }\mathbf{T}_{y}f)(x,\xi )=e^{2\pi i\omega  x}\mathbf{Z}f(x-y,\xi -\omega ),  \label{eq:zak-tf-shift}
\end{equation} \vspace{-0.3cm}

\noindent \emph{Unitarity:} for $f_{1},f_{2}\in \mathrm{L}^{2}(\mathbb{R})$ it holds 
\begin{equation}
\int_{0}^{1}\int_{0}^{1 }\mathbf{Z}f_{1}(x,\xi )\overline{
\mathbf{Z}f_2(x,\xi )}dxd\xi =\langle f_{1},f_{2}\rangle_{\mathrm{L}^{2}} .
\label{eq:unitary}
\end{equation}
\medskip

The following lemma shows that $\mathbf{V}_{g}$ can be easily computed as a linear combination of shifted Zak transforms of the window $g$.

\begin{lemma}\label{lem:stft-explicit}
Let $f,g\in {S}_0(\mathbb{R})$ and $\varphi\in S_N$. Then
\begin{equation}
\mathbf{V}_{g}(\mathbf{\Sigma}_{N}f)(x,\xi )=\sum_{n=0}^{N-1} \mathbf{P}f\left( \frac{n}{N}
\right) e^{-2\pi i\xi \frac{n}{N}}\mathbf{Z}\overline{g}\left( \frac{n}{N}-x,\xi
\right) ,  \label{eq:V_gf-double-per}
\end{equation}
and for $\varphi=\sum_{n=1}^Na_n\epsilon_n$
\begin{equation}
\mathbf{V}_{g}\varphi (x,\xi )=\sum_{n=0}^{N-1}a_{n}e^{-2\pi i\xi \frac{n}{N}}\mathbf{Z}
\overline{g}\left( \frac{n}{N}-x,\xi \right) .  \label{eq:V_gvarphi}
\end{equation}
\end{lemma}
\proof First, let us compute the STFT of the basis functions $\epsilon _{n}$
\begin{align}
\mathbf{V}_{g}\epsilon _{n}(x,\xi )&=\left\langle \epsilon _{n},\mathbf{M}_{\xi
}\mathbf{T}_{x}g\right\rangle _{\mathcal{S}_{0}^{\prime } \times \mathcal{S}_0}  =\sum_{k\in \mathbb{Z}}\left\langle \delta _{\frac{n}{N}+k},\mathbf{M}_{\xi
}\mathbf{T}_{x}g\right\rangle _{\mathcal{S}_{0}^{\prime } \times \mathcal{S}_0}  \notag \\
&=\sum_{k\in \mathbb{Z}}\overline{g\left( \frac{n}{N}+k-x\right) }e^{-2\pi
i\xi (\frac{n}{N}+k)} =e^{-2\pi i\xi \frac{n}{N}}\mathbf{Z}\overline{g}\left( \frac{n}{N}-x,\xi
\right) ,  \label{Basis}
\end{align}
For general $\varphi =\sum_{n=0}^{N-1}a_{n}\epsilon _{n}\in S_{N}$ one thus
gets \eqref{eq:V_gvarphi}.
Applying $\mathbf{V}_{g}$ to (\ref{eq:a_n}) from\ Lemma~\ref{lem:double-period-1} gives 
\begin{equation*}
\mathbf{V}_{g}(\Sigma _{N}f)(x,\xi )=\frac{1}{N}\sum_{n=0}^{N-1} \mathbf{P}f\left( 
\frac{n}{N}\right) \mathbf{V}_{g}\epsilon _{n}(x,\xi )
\end{equation*}
which combined with \eqref{Basis} yields \eqref{eq:V_gf-double-per}.
\pbox

\noindent With these basic observations, it is now straightforward to show Theorem~\ref{thm:finite-gabor}. For convenience, we repeat the statement here.\medskip

\newcounter{orgthm}
\setcounter{orgthm}{\value{theorem}}
\setcounter{theorem}{0}

\begin{theorem}  Let $f,g\in S_{0}(\mathbb{R})$ and let $\mathsf{f_N} = \mathbf{P}_N f,\mathsf{g_N} = \mathbf{P}_N g$.
Then, for $l,k\in 0,....,N-1$,
\begin{equation*}
\mathbf{V}_{g}( \mathbf{\Sigma} _{N}f)\left( \frac{k}{N},l\right) =N^{-1}\cdot\V_{\mathsf{g_N}}\mathsf{f_N}[k,l].
\end{equation*}
\end{theorem}
\proof
Setting $\xi =l\in 0,....,N-1$ and $x=\frac{k}{N},\ k\in 0,\ldots ,N-1,$
yields 
\begin{align*}
\mathbf{V}_{g}(\mathbf{\Sigma} _{N}f)\left(\frac{k}{N},l\right) &=\frac{1}{N}\sum_{n=0}^{N-1} \mathbf{P}
f\left( \frac{n}{N}\right) e^{-2\pi i\frac{nl}{N}}\mathbf{Z}\overline{g}\left( 
\frac{n-k}{N},l\right)  \\
&=\frac{1}{N}\sum_{n=0}^{N-1} \mathbf{P}f\left( \frac{n}{N}\right) e^{-2\pi i\frac{nl}{N}}\overline{ \mathbf{P}g}\left( \frac{n-k}{N}\right) 
\end{align*}
This shows that by definition  of the \dgt 
 we get the equality 
\begin{equation*}
\mathbf{V}_{g}(\mathbf{\Sigma }_{N}f)\left( \frac{k}{N},l\right) =\frac{1}{N}\V_{\mathsf{g_N}}\mathsf{f_N}[k,l].
\end{equation*}
\pbox

\setcounter{theorem}{\value{orgthm}}

{
\begin{remark}
In fact, we can express $\mathbf{V}_{g}\varphi (x,\xi )$, with $\varphi=\sum_{n=1}^Na_n\epsilon_n$, in terms of the \dgt for arbitrary $(x,\xi)\in\mathbb{T}^2_N$. Since the expression in Eq. \eqref{eq:V_gvarphi} relies on the $\bd L^2(\mathbb{R})$-Zak transform, doing so may be more convenient for implementations. Specifically, with $\mathrm{f_N}=(a_1,\ldots,a_N)^T\in \mathbb{C}^N$, we have 
\begin{equation}
\mathbf{V}_{g}\varphi (x,\xi )= \exp(-2\pi i r_\xi n_x/N)\cdot \mathbf{V}_{\mathrm{g}^{(r_x,r_\xi)}_\mathrm{N}}\mathrm{f_N}[n_x,m_\xi],\label{eq:V_pure_disc}
\end{equation}
where $x = n_x/N+r_x$ and, $\xi = m_\xi+r_\xi$, with $r_x\in [0,1/N)$, $r_\xi\in [0,1)$, and 
\[
\mathrm{g}_{\mathrm{N}}^{(r_x,r_\xi)} = \mathbf{P}_N (\bd M_{r_\xi} \bd T_{r_x} g).
\]
To show this, begin by inserting the definition of $\mathbf{Z}\overline{g}$ in Eq. \eqref{eq:V_gvarphi} to obtain
\[
 \begin{split}
   \mathbf{V}_{g}\varphi (x,\xi ) & = \sum_{n=0}^{N-1}a_{n}e^{-2\pi i\xi \frac{n}{N}}\cdot \sum_{k\in\mathbb{Z}} \overline{g}(n/N-x-k)e^{2\pi ik\xi}\\
    & = \sum_{n=0}^{N-1}a_{n}e^{-2\pi im_\xi \frac{n}{N}}\cdot \sum_{k\in\mathbb{Z}} \overline{g}\left(\frac{n-n_x}{N}-k-r_x\right)e^{-2\pi i(\frac{n}{N}-k)r_\xi}.
 \end{split}
\]
Now, note that
\[
\begin{split}
   \lefteqn{\sum_{k\in\mathbb{Z}} \overline{g}\left(\frac{n-n_x}{N}-k-r_x\right)e^{-2\pi i(\frac{n}{N}-k)r_\xi}}\\
   & = e^{-2\pi ir_\xi \frac{n_x}{N}}\cdot \overline{\left(\sum_{k\in\mathbb{Z}} g(\, \bullet\, -k-r_x) e^{2\pi i(\, \bullet\,-k)r_\xi}\right)}\left(\frac{n-n_x}{N}\right)\\
   & = e^{-2\pi ir_\xi \frac{n_x}{N}}\cdot \overline{\mathbf{P}_N (\bd M_{r_\xi} \bd T_{r_x} g)}[n-n_x] = e^{-2\pi ir_\xi \frac{n_x}{N}}\cdot \overline{\mathrm{g}_\mathrm{N}^{(r_x,r_\xi)}}[n-n_x].
 \end{split} 
\]
Together, we have
\[
  \mathbf{V}_{g}\varphi (x,\xi ) = e^{-2\pi ir_\xi \frac{n_x}{N}} \sum_{n=0}^{N-1}a_{n}e^{-2\pi im_\xi \frac{n}{N}}\cdot \overline{\mathrm{g}_\mathrm{N}^{(r_x,r_\xi)}}[n-n_x],
\]
proving \eqref{eq:V_pure_disc}, since $\mathrm{f_N}[n] = a_n$.
\end{remark}}

\subsection{Moyal-type and inversion formulas for the STFT on flat tori}

One of the most fundamental identities in time-frequency analysis is Moyal's formula ensuring that the STFT is an isometry from $\mathrm{L}^{2}(\R)$ to $\mathrm{L}^{2}(\R^{2})$
$$
\int_{\R^2}\mathbf{V}_{g_1}f_1(x,\xi)\overline{\mathbf{V}_{g_2}f_2(x,\xi)}dxd\xi= \langle f_1,f_2\rangle_{\mathrm{L}^{2}} \langle g_2,g_1\rangle_{\mathrm{L}^{2}},\quad f_1,f_2,g_1,g_2\in \mathrm{L}^{2}(\R).
$$
We will now show the toric equivalent of this identity ensuring that the STFT on the flat tori is a multiple of an isometry.
\begin{theorem}
\label{thm:moyal} If $\varphi_1,\varphi_2\in S_N$, and $g_1,g_2\in {S}
_0(\mathbb{R})$, then 
\begin{equation}\label{eq:moyal}
\int_{\mathbb{T}_N^2}\mathbf{V}_{g_1}\varphi_1(x,\xi)\overline{\mathbf{V}_{g_2}\varphi_2(x,\xi)}
dxd\xi=N\langle \varphi_1,\varphi_2\rangle_{S_N}\langle g_2,g_1\rangle_{\mathrm{L}^{2}}.
\end{equation}
\end{theorem}

\proof
For $\varphi _{1}=\sum_{n=0}^{N-1}a_{n}\epsilon _{n}$, and $\varphi
_{2}=\sum_{n=0}^{N-1}b_{n}\epsilon _{n}$, we  use Lemma~\ref{lem:stft-explicit} and the periodicity of the Zak transform in the frequency variable (\ref{eq:quasiperiodic}) to obtain
\begin{align*}
\int_{\mathbb{T}_{N}^2}\mathbf{V}_{g_{1}}\varphi _{1}&(x,\xi ) \overline{\mathbf{V}_{g_{2}}\varphi
_{2}(x,\xi )}dxd\xi
\\
&= \int_{0}^{N}\int_{0}^{1}\sum_{n,k=0}^{N-1}a_{n}\overline{b_{k}}e^{-2\pi
i\xi \frac{n-k}{N}}\mathbf{Z}\overline{g_{1}}\left( \frac{n}{N}-x,\xi \right) 
\overline{\mathbf{Z}\overline{g_{2}}\left( \frac{k}{N}-x,\xi \right) }dxd\xi  
\\
&= \sum_{l=0}^{N-1}\int_{0}^{1}\int_{0}^{1}\sum_{n,k=0}^{N-1}a_{n}\overline{b_{k}}e^{-2\pi
i(\xi +l)\frac{n-k}{N}}\cdot
\\
& \hspace{95pt}\mathbf{Z}\overline{g_{1}}\left( \frac{n}{N}-x,\xi
+l\right)\overline{\mathbf{Z}\overline{g_{2}}\left( \frac{k}{N}-x,\xi +l\right) }
dxd\xi  
\\
&= \int_{0}^{1}\int_{0}^{1}\sum_{n,k=0}^{N-1}a_{n}\overline{b_{k}}\left(
\sum_{l=0}^{N-1}e^{-2\pi il\frac{n-k}{N}}\right) e^{-2\pi i\xi\frac{n-k}{N}
}\cdot
\\
& \hspace{95pt}\mathbf{Z}\overline{g_{1}}\left( \frac{n}{N}-x,\xi \right) \overline{\mathbf{Z}
\overline{g_{2}}\left( \frac{k}{N}-x,\xi \right) }dxd\xi 
\\
&= N\sum_{n=0}^{N-1}a_{n}\overline{b_{n}}\int_{0}^{1}\int_{0}^{1}\mathbf{Z}
\overline{g_{1}}\left( \frac{n}{N}-x,\xi \right) \overline{\mathbf{Z}\overline{
g_{2}}\left( \frac{n}{N}-x,\xi \right) }dxd\xi ,
\end{align*}
where we used the basic fact $
\sum_{l=0}^{N-1}e^{-2\pi il\frac{k}{N}}=N\delta _{k,0}$, for $k = 0, \dots, N-1$, to obtain the final equality. Applying
consecutively (\ref{eq:zak-tf-shift}), (\ref{eq:quasiperiodic}) and (\ref
{eq:unitary}) to the integral above yields  
\begin{align*}
 \int_{0}^{1}\int_{0}^{1}&\mathbf{Z}
\overline{g_{1}}\left( \frac{n}{N}-x,\xi \right) \overline{\mathbf{Z}\overline{
g_{2}}\left( \frac{n}{N}-x,\xi \right) }dxd\xi
\\
&= \int_{0}^{1}\int_{0}^{1}\mathbf{Z}(T_{-
\frac{n}{N}}\overline{g_{1}})\left( -x,\xi \right) \overline{\mathbf{Z}(T_{-\frac{
n}{N}}\overline{g_{2}})\left( -x,\xi \right) }dxd\xi 
 \\
&= \int_{0}^{1}\int_{0}^{1}e^{-2\pi
i\xi }\mathbf{Z}(\mathbf{T}_{-\frac{n}{N}}\overline{g_{1}})\left( 1-x,\xi \right) 
\overline{e^{-2\pi i\xi }\mathbf{Z}(\mathbf{T}_{-\frac{n}{N}}\overline{g_{2}})\left(
1-x,\xi \right) }dxd\xi  
\\
&=\int_{0}^{1}\int_{0}^{1}\mathbf{Z}(\mathbf{T}_{-
\frac{n}{N}}\overline{g_{1}})\left( x,\xi \right) \overline{\mathbf{Z}(\mathbf{T}_{-\frac{n
}{N}}\overline{g_{2}})\left( x,\xi \right) }dxd\xi  
\\
&=\langle \mathbf{T}_{-\frac{n}{N}}\overline{
g_{1}},\mathbf{T}_{-\frac{n}{N}}\overline{g_{2}}\rangle_{\mathrm{L}^{2}} = \langle g_{2},g_{1}\rangle_{\mathrm{L}^{2}} ,
\end{align*}
which concludes the proof.
\pbox

\begin{remark}
Using \eqref{eq:sigma_N-f,g} we thus have shown that $\big\{\mathbf{\Sigma}_N (\pi(x,\xi)g)\big\}_{(x,\xi)\in\mathbb{T}_N^2}$ is a continuous tight frame for $S_N$, with bound $N \cdot \Vert g\Vert _{\mathrm{L}^{2}} $. 
\end{remark}

\noindent From Theorem~\ref{thm:moyal} and \eqref{eq:sigma_N-f,g} we can now
derive two inversion formulas.

\begin{theorem}
\label{thm:inversion} Let $g\in {S}_0(\mathbb{R})\backslash\{0\}$. For every $
\varphi\in S_N$ and every $f\in {S}_0(\mathbb{R})$, it holds 
\begin{equation}  \label{eq:inversion-SN}
\varphi=\frac{1}{N\|g\|^2_2}\sum_{n=0}^{N-1}\left(\int_{\mathbb{T}
_N^2}\mathbf{V}_g\varphi (x,\xi)e^{2\pi i\xi\frac{n}{N}} \mathbf{Z}g\left(\frac{n}{N}
-x,-\xi\right)dxd\xi \right)\epsilon_n ,
\end{equation}
and 
\begin{equation}  \label{eq:inversion-db-period}
\mathbf{\Sigma} _{N}f=\frac{1}{\|g\|^2_2 }\int_{\mathbb{T}_N^2}\mathbf{V}_{g}(\mathbf{\Sigma} _{N}f)(x,\xi
)\mathbf{\Sigma}_{N}\big(\pi (x,\xi )g\big)dxd\xi,
\end{equation}
where the integral is understood in the weak-$\ast$ sense. 
\end{theorem}

\proof
Let $\varphi=\sum_{n=0}^{N-1}a_n\epsilon_n$. Then by Theorem~\ref{thm:moyal}
we know 
\begin{align*}
a_n&=\langle \varphi,\epsilon_n\rangle_{S_N}=\frac{1}{N\|g\|^2_2}\int_{
\mathbb{T}_N^2}\mathbf{V}_g\varphi(x,\xi)\overline{\mathbf{V}_g\epsilon_n(x,\xi)}dxd\xi 
\\
&= \frac{1}{N\|g\|^2_2}\int_{\mathbb{T}_N^2}\mathbf{V}_g\varphi(x,\xi)e^{2\pi i\xi\frac{n
}{N}}\overline{\mathbf{Z}\overline{g}\left(\frac{n}{N}-x,\xi\right)}dxd\xi 
\\
&= \frac{1}{N\|g\|^2_2}\int_{\mathbb{T}_N^2}\mathbf{V}_g\varphi(x,\xi)e^{2\pi i\xi\frac{n
}{N}}\mathbf{Z}g\left(\frac{n}{N}-x,-\xi\right)dxd\xi,
\end{align*}
where we used (\ref{eq:V_gvarphi}) and the fact that $\overline{\mathbf{Z}
g(x,\xi)}=\mathbf{Z}\overline{g}(x,-\xi)$. To show the second identity we
first observe that $\mathbf{\Sigma}_N$ satisfies 
\begin{equation*}
\langle \mathbf{\Sigma}_N f,g\rangle_{ {S}_0^\prime \times  {S
}_0 }=\langle f,\mathbf{\Sigma}_N g\rangle_{ {S}_0 \times  {S}_0^\prime },\quad f,g\in S_0(\R).
\end{equation*}
Let $f,h\in\mathcal{S}_0( \mathbb{R})$. The result thus follows from (\ref
{eq:sigma_N-f,g}) and \eqref{eq:moyal} as 
\begin{align*}
\langle \mathbf{\Sigma}_N f,h\rangle_{ {S}_0^\prime \times  {S}_0 }&=\frac{1}{\|g\|^2_2}\int_{\mathbb{T_N}^2}\mathbf{V}_g(\mathbf{\Sigma}_N f)(x,\xi)
\langle \pi(x,\xi) g, \mathbf{\Sigma}_N h\rangle_{ {S}_0 \times 
 {S}_0^\prime }dxd\xi \\
&=\frac{1}{\|g\|^2_2}\int_{\mathbb{T_N}^2}\mathbf{V}_g(\mathbf{\Sigma}_N f)(x,\xi) \langle \mathbf{\Sigma}_N
\big(\pi(x,\xi) g\big), h\rangle_{ {S}_0^\prime \times 
 {S}_0 }dxd\xi.
\end{align*}
\hfill $\Box$

\begin{remark}
Note that coorbit theory \cite{fe06,FeGr1} also guarantees an inversion of the short-time Fourier transform on $ {S}_0^\prime(\R)$ as $\mathbf{V}_g^\ast \mathbf{V}_g=\|g\|^2 \mathbf{I}_{ {S}_0^\prime}$. The benefit of our point of view is however that the coefficients in \eqref{eq:inversion-SN} can directly be calculated without resorting to weakly defined integrals.
\end{remark}

\subsection{Reproducing kernel}

We now prove that the range of the STFT restricted to $S_N$ is an $N$-dimensional reproducing kernel Hilbert space (RKHS) of $\mathrm{L}^{2}(\mathbb{T}_N^2)$ and give several expressions for its reproducing kernel. 
\begin{proposition}
The space $\mathbf{V}_{g}(S_N)\subset \mathrm{L}^{2}(\mathbb{T}_N^2)$ is  a RKHS and its kernel is given by 
\begin{align}
K_{g }\big((x^{\prime },\xi ^{\prime }),(x,\xi )\big)  
& =\frac{1
}{N\Vert g\Vert ^{2}_2}\sum_{n=0}^{N-1}e^{2\pi i(\xi -\xi ^{\prime })\frac{n}{N
}}\overline{\mathbf{Z}\overline{g}\left( \frac{n}{N}-x,\xi \right) }\mathbf{Z}
\overline{g}\left( \frac{n}{N}-x^{\prime },\xi ^{\prime }\right) 
\notag
\\
& =\frac{N}{\Vert g\Vert ^{2}_2}\langle\mathbf{ \Sigma }_{N}\big(\pi (x,\xi )g\big)
,\mathbf{\Sigma} _{N}\big(\pi (x^{\prime },\xi ^{\prime })g\big)\rangle _{S_{N}}
\label{eq:RK-1} 
\\
&=\frac{1}{\Vert g\Vert ^{2}_2}\sum_{k_{1},k_{2}\in \mathbb{Z}
} e^{-2\pi i k_1 \xi}\langle \pi (x+k_1,\xi+Nk_2 )g,\pi (x^{\prime
},\xi ^{\prime })g\rangle_{\mathrm{L}^{2}} .
\notag 
\end{align}
\end{proposition}
\proof
Using consecutively \eqref{eq:V_gvarphi}, Cauchy-Schwarz inequality, \cite[Lemma~8.2.1]{Charly} and  \eqref{eq:moyal} yields
\begin{align*}
|\mathbf{V}_g\varphi(x,\xi)|&\leq \|\varphi\|_{S_N}\left(\sum_{n=0}^{N-1} \left|\mathbf{Z}\overline{g}\left(\frac{n}{N}-x,\xi\right)\right|^2\right)^{1/2}
\\
&\leq \|\varphi\|_{S_N}\sqrt{N}\|\mathbf{Z}\overline{g}\|_\infty=\frac{\|\mathbf{Z}\overline{g}\|_\infty}{\|g\|_2}\|\mathbf{V}_g\varphi\|_{\mathrm{L}^{2}(\mathbb{T}_N)},
\end{align*}
meaning that point evaluation is continuous and $V_g(S_N)$ is a RKHS.

By Theorem~\ref{thm:inversion} and Lemma~\ref{lem:stft-explicit}, it follows 
\begin{align*}
\mathbf{V}_g&\varphi(x^\prime,\xi^\prime)=\hspace{-0.05cm}\frac{1}{N\|g\|^2_2}\sum_{n=0}^{N-1}\int_{
\mathbb{T}_N^2}\hspace{-0.1cm}\mathbf{V}_g\varphi(x,\xi)e^{2\pi i \xi\frac{n}{N}}\mathbf{Z}g\left(\frac{n}{N}
-x,-\xi\right)\hspace{-0.05cm}\mathbf{V}_g\epsilon_n(x^\prime,\xi^\prime)dxd\xi 
\\
&=\hspace{-0.05cm}\int_{\mathbb{T}_N^2}\hspace{-0.1cm}\mathbf{V}_g\varphi(x,\xi)\frac{1}{N\|g\|_2^2}\sum_{n=0}^{N-1}e^{2
\pi i (\xi-\xi^\prime)\frac{n}{N}}\overline{\mathbf{Z}\overline{g}\left(\frac{n}{N}
-x,\xi\right)}
\mathbf{Z}\overline{g}\left(\frac{n}{N}-x^\prime,\xi^\prime\right)dxd
\xi,
\end{align*}
as well as by \eqref{eq:inversion-db-period} and \eqref{eq:sigma_N-f,g}
\begin{align*}
\mathbf{V}_g(\mathbf{\Sigma}_Nf)&(x^\prime,\xi^\prime)\\ &=\frac{1}{\|g\|^2_2}\int_{\mathbb{T}
_N^2}\mathbf{V}_g\mathbf{\Sigma}_Nf(x,\xi)\langle \mathbf{\Sigma}_N\big(\pi(x,\xi)g\big)
,\pi(x^\prime,\xi^\prime)g\rangle_{\mathcal{S}_0^\prime(\mathbb{R})\times
\mathcal{S}(\mathbb{R})}dxd\xi \\
&=\frac{N}{\|g\|^2_2}\int_{\mathbb{T}_N^2}\mathbf{V}_g\mathbf{\Sigma}_Nf(x,\xi)\langle \mathbf{\Sigma}_N
\big(\pi(x,\xi)g\big), \mathbf{\Sigma}_N\big(\pi(x^\prime,\xi^\prime)g\big)
\rangle_{S_N}dxd\xi.
\end{align*}
Therefore, the first two identities of \eqref{eq:RK-1} hold.
Moreover, 
\begin{align*}
K_{g}\left( (x^{\prime },\xi ^{\prime }),(x,\xi )\right)  &=
\frac{N}{ \Vert g\Vert ^{2}_2}\langle\mathbf{\Sigma} _{N}(\pi (x ,\xi )g ) ,\mathbf{\Sigma}_N(\pi(x^\prime,\xi^\prime)g)\rangle_{S_N}
\\
&=\frac{1}{\Vert g\Vert ^{2}_2}\langle\mathbf{ \Sigma}_N (\pi(x ,\xi  )g),\pi(x^\prime,\xi^\prime)g\rangle_{\mathcal{S}_0^\prime\times\mathcal{S}_0}
\\
&=\frac{1}{\Vert g\Vert ^{2}_2}\sum_{k_{1},k_{2}\in \mathbb{Z}%
}e^{-2\pi ik_{1}\xi }\langle\pi(x +k _{1},\xi 
+Nk_{2})g,\mathbf{\pi}(x^\prime,\xi^\prime)g\rangle_{\mathrm{L}^{2}} .
\end{align*}
\pbox

\subsection{Examples of the STFT on flat tori}

In this section, we consider explicit calculations of the objects discussed in the previous sections using the  non-normalized dilated Gaussian
\begin{equation*}
h_{0}^\lambda(t):=e^{-\pi \lambda t^{2}},\quad \lambda>0.
\end{equation*}
The reason we introduce the additional dilation is that, having  the connection to finite Gabor systems (Theorem~\ref{thm:finite-gabor}) in mind, we would like to impose a certain degree of  localization of $\mathbf{P}h_0^\lambda$. This is only guaranteed if $\lambda>0$ is chosen large enough, see Figure~\ref{fig:1}.

Following \eqref{eq:V_gvarphi}, the STFT of a basis function $\epsilon_n$ is given by
\begin{align}
\mathbf{V}_{h_{0}^\lambda}\epsilon _{n}(x,\xi ) &=\,e^{-2\pi i\frac{n}{N}\xi }\mathbf{Z}\overline{
h_{0}^\lambda}\left( \frac{n}{N}-x,\xi \right)   \notag  \label{eq:h0-basis-Vg} \\
&= e^{-2\pi i\frac{n}{N}\xi }\sum_{k\in \mathbb{Z}}e^{-\pi \lambda (x-\frac{n
}{N}+k)^{2}}e^{2\pi i\xi k}  \notag \\
&= e^{-\pi \lambda(x-\frac{n}{N})^{2}-2\pi i\xi \frac{n}{N}}\sum_{k\in 
\mathbb{Z}}e^{-\pi  \lambda k^{2}}e^{2\pi ik\left( \xi +i\lambda(x-\frac{n}{N})\right) } 
\notag \\
&= e^{-\pi \lambda x^{2}}e^{-\pi \lambda \left(\frac{n}{N}\right) ^{2}}e^{2\pi \overline{z_\lambda}\frac{n}{N}}\vartheta \left( i\left( \overline{z_\lambda}-\lambda\frac{n}{N}
\right) ,i\lambda\right) ,
\end{align}
where $z_\lambda:=\lambda x+i\xi $ and the Jacobi theta function $\vartheta $ is defined as 
\begin{equation*}
\vartheta (z,\tau ):=\sum_{k\in \mathbb{Z}}e^{\pi ik^{2}\tau +2\pi ikz},\quad \text{Im}(\tau)>0.
\end{equation*}

\begin{center}
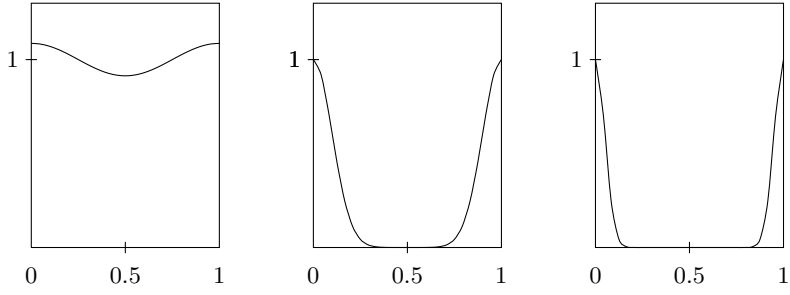
\begin{figure}
\begin{tikzpicture}[scale=2.5]
\draw (0,0)--(0,1.3)--(1,1.3)--(1,0)--(0,0);
\draw (-0.03,1)--(0.03,1);
\node at (1.4,1) {\small{1}};
\node at (-0.1,1) {\small{1}};
\draw[domain=0:1, smooth, variable=\x, black] plot ({\x}, {exp(-pi*\x^2)+exp(-pi*(\x-1)^2)+exp(-pi*(\x-2)^2)+exp(-pi*(\x+1)^2)+exp(-pi*(\x+2)^2});

\draw (1.5,0)--(1.5,1.3)--(2.5,1.3)--(2.5,0)--(1.5,0);
\draw (1.47,1)--(1.53,1);
\node at (1.4,1) {\small{1}};
\draw[domain=1.5:2.5, smooth, variable=\x, black] plot ({\x}, {exp(-pi*16*(\x-1.5)^2)+exp(-pi*16*(\x-2.5)^2)+exp(-pi*16*(\x-3.5)^2)+exp(-pi*16*(\x-0.5)^2)+exp(-pi*16*(\x+0.5)^2});
\node at (0,-0.15) {\small{0}};
\node at (2.5,-0.15) {\small{1}};
\node at (1,-0.15) {\small{1}};
\node at (1.5,-0.15) {\small{0}};

\draw (3,0)--(3,1.3)--(4,1.3)--(4,0)--(3,0);
\draw (2.97,1)--(3.03,1);
\node at (2.9,1) {\small{1}};
\draw[domain=3:4, smooth, variable=\x, black] plot ({\x}, {exp(-pi*64*(\x-3)^2)+exp(-pi*64*(\x-4)^2)+exp(-pi*64*(\x-5)^2)+exp(-pi*64*(\x-2)^2)+exp(-pi*64*(\x-1)^2});
\node at (3,-0.15) {\small{0}};
\node at (4,-0.15) {\small{1}};
\node at (0.5,-0.15) {\small{0.5}};
\draw (0.5,-0.03)--(0.5,0.03);
\node at (2,-0.15) {\small{0.5}};
\draw (2,-0.03)--(2,0.03);
\node at (3.5,-0.15) {\small{0.5}};
\draw (3.5,-0.03)--(3.5,0.03);
\end{tikzpicture}
\caption{$\mathbf{P}h_0^\lambda$, for $\lambda=1$, (left) $\lambda=16$ (middle), and $\lambda=64$ (right). }\label{fig:1}
\end{figure}
\end{center}

\noindent Moreover, we explicitly calculate the reproducing kernel  $K_{h_{0}^\lambda,\mathbb{T}^2_{N}}$. For this purpose, the third equation of \eqref{eq:RK-1} is the most convenient representation of the kernel. Let us write $w_\lambda=\lambda x'+i\xi'$. By \cite[Lemma~1.5.2]{Charly} 
\begin{align*}\label{eq:gauss-tf}
\langle \pi(x,\xi )h_{0}^\lambda,\pi(x^\prime,\xi^\prime )h_0^\lambda\rangle  
 =(2\lambda)^{-\frac{1}{2}}e^{\pi i(\xi -\xi^\prime)(x+x^\prime)}e^{-\frac{\pi }{2}\left[
\lambda ( x-x^\prime ) ^{2}+\frac{1}{\lambda}(\xi -\xi^\prime)^{2}\right] }.
\end{align*}
from which we deduce 
\begin{align*}
K&_{h_{0}^\lambda}(w_\lambda,z_\lambda) =\frac{1}{\|h_0^\lambda\|^2}\sum_{k_1,k_2\in\mathbb{Z}}e^{-2\pi ik_1\xi} \langle \pi(x+k_1,\xi+Nk_2)h_0^\lambda,\pi(x^\prime,\xi^\prime)h_0^\lambda\rangle
\\
&=  e^{\pi i(\xi -\xi^\prime)(x+x')}\hspace{-0.1cm}\sum_{k_{1},k_{2}\in \mathbb{Z}}e^{\frac{\pi }{2}  \left[-2i(\xi +\xi^\prime)k_1+ 2i(x+x^\prime)Nk_2-\lambda ( x-x^\prime +k_1 ) ^{2}-\frac{1}{\lambda}(\xi -\xi^\prime+Nk_2)^{2}\right] } 
\\
&= e^{\pi\left[i(\xi -\xi^\prime)(x+x')-\frac{\lambda}{2}(x-x')^2-\frac{1}{2\lambda}(\xi-\xi')^2\right]}\hspace{-0.1cm}\sum_{k_{1} \in \mathbb{Z}}e^{ -\pi   k_1 [\lambda(x-x')+i(\xi+\xi')] -\pi \frac{ \lambda k_1^2}{2}}  
\\
&\hspace{2cm}
\sum_{k_{2} \in \mathbb{Z}} e^{ {\pi}i\frac{  Nk_2}{\lambda} \left[  \lambda (x+x^\prime)  +i (\xi -\xi^\prime)\right]-\pi\frac{N^2k_2^{2}}{2\lambda} } 
\\
&= e^{- \pi \lambda  (x^2+ (x^\prime)^2) }e^{\frac{\pi }{2\lambda }(z_\lambda+\overline{w_\lambda})^2}\vartheta \left( i\frac{z_\lambda- \overline{w_\lambda}}{2}, \frac{i\lambda }{
2}\right) \vartheta \left( \frac{(z_\lambda+\overline{w_\lambda})N}{2\lambda},\frac{iN^{2}}{2\lambda}\right).
\end{align*}

\section{A Bargmann-type transform and finite Gabor Frames}

Given $\varphi =\sum_{n=0}^{N-1}a_{n}\epsilon _{n}\in S_{N}$\ we can define
a Bargmann-type transform  by $ {\mathbf{B}}_{(\lambda,N)}\varphi
(z)=\mathbf{V}_{h_{0}^\lambda}\varphi \left(\frac{x}{\lambda},-\xi\right)e^{\pi  x^{2}/\lambda}$, $z=x+i\xi$, or equivalently using \eqref{eq:h0-basis-Vg}  
\begin{equation*}
\mathbf{B}_{(\lambda,N)}\varphi (z)= \sum_{n=0}^{N-1}a_{n}e^{-\pi \lambda\left( \frac{n
}{N}\right) ^{2}}e^{2\pi  z\frac{n}{N}}\vartheta \left( i\left( z-\lambda\frac{n}{N}
\right) ,i\lambda\right) \text{.}
\end{equation*}
By the quasi-periodicity of $\mathbf{V}_{h_{0}^\lambda}\varphi $  \eqref{eq:quasiper}, we then get  that
\begin{equation}
\mathbf{B}_{(\lambda,N)}\varphi (z+\lambda n)=e^{\pi \lambda n^{2}+2\pi zn}\mathbf{B}_{(\lambda,N)}\varphi
(z),\quad n\in \mathbb{Z},  \label{eq:barg-time-shift}
\end{equation}
as well as  
\begin{equation}
\mathbf{B}_{(\lambda,N)}\varphi (z+iNm)=\mathbf{B}_{(\lambda,N)}\varphi (z),\quad m\in \mathbb{
Z}.  \label{eq:barg-freq-shift}
\end{equation}
This can also be derived using the quasi-periodicity of the Jacobi theta
function 
\begin{equation}\label{eq:jacobi-quasi}
\vartheta (z+n+\tau m,\tau )=e^{-\pi i\tau m^{2}}e^{-2\pi imz}\vartheta
(z,\tau ),\quad n,m\in \mathbb{Z}.
\end{equation}
We will now show that the range of the Bargmann-type transform ${\mathbf{B}}_{(\lambda,N)}$ is precisely the space of analytic functions satisfying \eqref{eq:barg-time-shift} and \eqref{eq:barg-freq-shift}. This follows from the following result.

\begin{lemma}
\label{dimension}
The space of entire functions satisfying the periodicity conditions \eqref{eq:barg-time-shift} and \eqref{eq:barg-freq-shift} is $N$-dimensional. Thus, any such function $F$ can be written as the linear
combination of $N$ orthogonal functions
\begin{equation*}
F(z)=\sum_{n=0}^{N-1}a_{n}e^{-\pi \lambda \left( \frac{n}{N}\right)
^{2}}e^{2\pi z\frac{n}{N}}\vartheta \left( i\left( z-\lambda \frac{n}{N}
\right) ,i\lambda \right) 
\end{equation*}
and $F(z)=\mathbf{B}_{(\lambda ,N)}\varphi (z)$, for some $\varphi
=\sum_{n=0}^{N-1}a_{n}\epsilon _{n}$.
\end{lemma}
\proof Let $F$ be analytic and satisfy \eqref{eq:barg-time-shift} and \eqref{eq:barg-freq-shift}. Since $F$ is $N$-periodic with respect to purely imaginary shifts, we can formally write $F$ as a Fourier series 
$$
F(z)=\sum_{k\in\mathbb{Z}} c_k e^{2\pi k z/N}.
$$
Plugging this expression into \eqref{eq:barg-time-shift} yields 
\begin{align*}
\sum_{k\in\mathbb{Z}} c_k e^{2\pi k (z+\lambda n)/N}&=e^{\pi \lambda  n^{2}+2\pi  zn}\sum_{k\in\mathbb{Z}} c_k e^{2\pi k z/N}=\sum_{k\in\mathbb{Z}} c_k e^{\pi \lambda  n^{2}} e^{2\pi (k+nN) z/N}
\\
&=\sum_{k\in\mathbb{Z}} c_{k-nN}e^{\pi\lambda n^2}e^{2\pi kz/N},
\end{align*}
which shows that the coefficients $c_k$ satisfy
\begin{equation}\label{eq:c_k}
c_{k+nN}=c_k e^{-\pi \lambda (n^2+2kn/N  )},\quad k\in\{0,...,N-1\},\ n\in\mathbb{Z}.
\end{equation}
There are therefore exactly $N$ coefficients $c_0,...,c_{N-1}$ that can be chosen freely and the other coefficients are determined by \eqref{eq:c_k}. Equation \eqref{eq:c_k} also shows that there exist constants $a,C>0$ such that  $|c_k|\leq C e^{-a(k/N)^2}$ which implies that $F$ coincides with its Fourier series, see, e.g., \cite[Exercise~I.4.4]{katz}. 

The orthogonality of the functions
\begin{equation*}
e^{-\pi \lambda \left( \frac{n}{N}\right) ^{2}}e^{2\pi z\frac{n}{N}
}\vartheta \left( i\left( z-\lambda \frac{n}{N}\right) ,i\lambda \right) =
\mathbf{V}_{h_{0}^{\lambda }}\epsilon _{n}\left( \frac{x}{\lambda },-\xi
\right) e^{\pi x^{2}/\lambda }
\end{equation*}
follows from $\left\langle \epsilon _{n},\epsilon _{m}\right\rangle
_{S_{N}}=\delta _{n,m}$ and Moyal's formula (\ref{eq:moyal}).
\pbox

\noindent The periodicity conditions \eqref{eq:barg-time-shift} and \eqref{eq:barg-freq-shift} will also allow us to count the number of zeros of 
$\mathbf{B}_{(\lambda,N)}\varphi $ in the torus $\mathbb{T}_{\lambda,N}^2:=\mathbb{R}^2/(\lambda\mathbb{Z}\times N\mathbb{Z})$, which coincides with the number or zeros of $\mathbf{V}_{h_0^\lambda}\varphi$ in $\mathbb{T}_{N}^2$, and to obtain a constraint these zeros must satisfy. Our arguments are inspired by \cite{LV}.

\begin{proposition}
\label{thm:N-zeros} Let $\varphi \in S_N\backslash\{0\}$. The function $\mathbf{B}_{(\lambda,N)}\varphi$ has exactly $N$ zeros (counted with their
multiplicities) on the torus $\mathbb{T}_{\lambda,N}^2$. Moreover, the $N$ zeros $z_1,...,z_N\in\mathbb{T}_{\lambda,N}^2$ satisfy
\begin{align}\label{eq:sum-zeros}
\sum_{k=1 }^Nz_k &=\frac{N\lambda }{2}+\lambda n+i\left(\frac{N^{2}}{2}+Nm\right),\quad \text{for some } n,m\in\mathbb{Z}.
\end{align}
\end{proposition}
\proof
Let the curve $\Gamma=\partial \mathbb{T}_{\lambda,N}^2$ be positively oriented and assume for now that $\Gamma$ contains no zero of $
\mathbf{B}_{(\lambda,N)}\varphi$. As $
\mathbf{B}_{(\lambda,N)}\varphi$ is an analytic function on $\mathbb{C}$, it follows  by
Cauchy's argument principle (see, e.g. \cite[Sec. 5.2]{ahl}) that the number of zeros of $\mathbf{B}_{(\lambda,N)}\varphi$ (counted with their multiplicities) is
given by 
\begin{align*}
\frac{1}{2\pi i }\int_\Gamma \frac{(\mathbf{B}_{(\lambda,N)}\varphi)^\prime}{\mathbf{B}_{(\lambda,N)}\varphi}d\gamma&=\frac{1}{2\pi i }\left[ \int_0^\lambda\frac{(\mathbf{B}
_{(\lambda,N)}\varphi)^\prime(t)}{\mathbf{B}_{(\lambda,N)}\varphi(t)}dt+i\int_0^N\frac{(\mathbf{B}
_{(\lambda,N)}\varphi)^\prime(\lambda+it)}{\mathbf{B}_{(\lambda,N)}\varphi(\lambda+it)}dt\right.  
\\
&\hspace{0.5cm}\left. -\int_0^\lambda\frac{(\mathbf{B}_{(\lambda,N)}\varphi)^\prime(t+iN)
}{\mathbf{B}_{(\lambda,N)}\varphi(t+iN)}dt-i\int_0^N\frac{(\mathbf{B}
_{(\lambda,N)}\varphi)^\prime(it)}{\mathbf{B}_{(\lambda,N)}\varphi(it)}dt \right] 
\\
&= \frac{1}{2\pi }\hspace{-0.05cm}\left[ \int_0^N\frac{(\mathbf{B}_{(\lambda,N)}\varphi)^\prime(\lambda+it)}{
\mathbf{B}_{(\lambda,N)}\varphi(\lambda+it)}dt-\hspace{-0.05cm}\int_0^N\frac{(\mathbf{B}_{(\lambda,N)}\varphi)^\prime(it)
}{\mathbf{B}_{(\lambda,N)}\varphi(it)}dt \right]\hspace{-0.05cm},
\end{align*}
where we used \eqref{eq:barg-freq-shift}. By \eqref{eq:barg-time-shift}
we have 
\begin{equation*}
(\mathbf{B}_{(\lambda,N)}\varphi)^\prime(z+\lambda)=2\pi e^{\pi \lambda +2\pi z }\mathbf{B}%
_{(\lambda,N)} \varphi(z)+e^{\pi \lambda +2\pi z }(\mathbf{B}_{(\lambda,N)}\varphi)^\prime(z),
\end{equation*}
and consequently 
\begin{align*}
\frac{1}{2\pi i } \int_\Gamma \frac{(\mathbf{B}_{(\lambda,N)}\varphi)^\prime}{\mathbf{B}_{(\lambda,N)}\varphi}d\gamma
 = \frac{1}{2\pi }&\left[ \int_0^N\frac{ e^{\pi\lambda +2\pi it }\big(2\pi
\mathbf{B}_{(\lambda,N)}\varphi(it)+ (\mathbf{B}_{(\lambda,N)}\varphi)^\prime(it)\big)}{ e^{\pi\lambda
+2\pi it }\mathbf{B}_{(\lambda,N)}\varphi(it)}dt\right.
\\
& \left.-\int_0^N\frac{(\mathbf{B}
_{(\lambda,N)}\varphi)^\prime(it)}{\mathbf{B}_{(\lambda,N)}\varphi(it)}dt \right]  =\frac{1}{2\pi}\int_0^N2\pi dt=N.
\end{align*}
Let $z_1,...,z_N$ be the $N$ zeros of $B_{(\lambda,N)}\varphi$. 
Then, using again \eqref{eq:barg-time-shift}, \eqref{eq:barg-freq-shift} and Cauchy's argument principle (see, e.g. \cite[Sec. 5.2, (49)]{ahl}) , we get with $\mathrm{id} \colon \mathbb{C} \rightarrow \mathbb{C}, z \mapsto z$
\begin{align*}
\sum_{k=1}^N z_k=& \frac{1}{2\pi i}\int_{\Gamma }{\mathrm{id}}\cdot \frac{(\mathbf{B}_{(\lambda ,N)}\varphi
)^{\prime }}{\mathbf{B}_{(\lambda ,N)}\varphi }d\gamma  
\\
 =&\frac{1}{2\pi i}\hspace{-0.05cm}\left[\int_{0}^{\lambda }\frac{t(\mathbf{B}_{(\lambda
,N)}\varphi )^{\prime }(t)}{\mathbf{B}_{(\lambda ,N)}\varphi (t)}
dt+i\int_{0}^{N}(\lambda +it)\frac{(\mathbf{B}_{(\lambda ,N)}\varphi
)^{\prime }(\lambda +it)}{\mathbf{B}_{(\lambda ,N)}\varphi (\lambda +it)}dt\right.
\\
& \hspace{0.8cm} \left.-\int_{0}^{\lambda }\frac{(t+iN)(\mathbf{B}_{(\lambda ,N)}\varphi
)^{\prime }(t+iN)}{\mathbf{B}_{(\lambda ,N)}\varphi (t+iN)}dt-i\int_{0}^{N}
\frac{it(\mathbf{B}_{(\lambda ,N)}\varphi )^{\prime }(it)}{\mathbf{B}
_{(\lambda ,N)}\varphi (it)}dt\right] 
\\
  =&\frac{1}{2\pi } \hspace{-0.05cm}\left[ -N\int_{0}^{\lambda }\frac{(\mathbf{
B}_{(\lambda ,N)}\varphi )^{\prime }(t)}{\mathbf{B}_{(\lambda ,N)}\varphi (t)
}dt- \int_{0}^{N}\frac{it(\mathbf{B}_{(\lambda ,N)}\varphi
)^{\prime }(it)}{\mathbf{B}_{(\lambda ,N)}\varphi (it)}dt\right.
\\
&\hspace{0.8cm}\left.+\int_{0}^{N}\frac{(\lambda +it)(\mathbf{B}_{(\lambda ,N)}\varphi
)^{\prime }(\lambda +it)}{\mathbf{B}_{(\lambda ,N)}\varphi (\lambda +it)}dt\right]  
\\
 =&\frac{1}{2\pi }\hspace{-0.05cm} \left[ -N\int_{0}^{\lambda }\frac{(\mathbf{
B}_{(\lambda ,N)}\varphi )^{\prime }(t)}{\mathbf{B}_{(\lambda ,N)}\varphi (t)
}dt -\hspace{-0.05cm} \int_{0}^{N} \frac{it(\mathbf{B}
_{(\lambda ,N)}\varphi )^{\prime }(it)}{\mathbf{B}_{(\lambda ,N)}\varphi (it)
}dt\right.
\\
&\hspace{0.8cm}\left.+\int_{0}^{N}(\lambda +it)\frac{2\pi \mathbf{B}_{(\lambda ,N)}\varphi
(it)+(\mathbf{B}_{(\lambda ,N)}\varphi )^{\prime }(it)}{\mathbf{B}_{(\lambda
,N)}\varphi (it)}dt\right] 
 \\
 =&\frac{1}{2\pi }\hspace{-0.05cm} \left[\lambda\hspace{-0.05cm} \int_{0}^{N}\hspace{-0.05cm}\frac{(\mathbf{B}_{(\lambda ,N)}\varphi )^{\prime
}(it)}{\mathbf{B}_{(\lambda ,N)}\varphi (it)}dt\hspace{-0.05cm}-\hspace{-0.05cm} N\hspace{-0.05cm}\int_{0}^{\lambda }\hspace{-0.05cm}\frac{(\mathbf{
B}_{(\lambda ,N)}\varphi )^{\prime }(t)}{\mathbf{B}_{(\lambda ,N)}\varphi (t)
}dt\hspace{-0.05cm}+\hspace{-0.05cm}2\pi\int_{0}^{N}(\lambda
+it)dt \right]\hspace{-0.05cm}. 
\end{align*}
Let us now think of  $t\mapsto \mathbf{B}_{(\lambda ,N)}\varphi (t)$ as a parametrization of  the curve   $\Gamma_1=\{\mathbf{B}_{(\lambda ,N)}\varphi (t)\}_{t\in[0,\lambda]}$, i.e. 
\[
\int_{0}^{\lambda }\frac{(\mathbf{B}_{(\lambda ,N)}\varphi )^{\prime }(t)}{
\mathbf{B}_{(\lambda ,N)}\varphi (t)}dt=\int_{\Gamma_1}\frac{1}{z}d\gamma.
\]
 Let $\Gamma_2=\big\{\mathbf{B}_{(\lambda ,N)}\varphi (0)(e^{\pi\lambda}(1-t) +t)\big\}_{t\in[0,1]}$ be the line segment that connects the endpoints of $\Gamma_1$.  Therefore, $\Gamma_0=\Gamma_1\cup\Gamma_2$ is a closed continuous curve. By the residue theorem, it follows that 
\[
\int_{\Gamma_0}\frac{1}{z}d\gamma=2\pi ik,
\]
for some  $k\in \mathbb{Z}$. Therefore 
\begin{align*}
    \int_{0}^{\lambda }\frac{(\mathbf{B}_{(\lambda ,N)}\varphi )^{\prime }(t)}{
\mathbf{B}_{(\lambda ,N)}\varphi (t)}dt&=-\int_{\Gamma_2}\frac{1}{z}d\gamma+2\pi ik=\int_0^1 \frac{e^{\pi\lambda}-1}{e^{\pi\lambda}-t(e^{\pi\lambda}-1)}dt+2\pi ik
\\
&=\int_0^{e^{\pi\lambda}-1} \frac{1}{e^{\pi\lambda}-t}dt+2\pi ik=\int_1^{e^{\pi\lambda}} \frac{1}{t}dt+2\pi ik=\pi \lambda +2\pi ik.
\end{align*}
Similarly,  $\Gamma_3=\{\mathbf{
B}_{(\lambda ,N)}\varphi (it)\}_{t\in[0,N]}$ is a closed and continuous curve. Let $n\in\mathbb{Z}$ be the winding number of $\Gamma_3$ around the origin. Then
\[
\int_{0}^{N}\frac{(\mathbf{B}_{(\lambda ,N)}\varphi )^{\prime }(it)}{\mathbf{
B}_{(\lambda ,N)}\varphi (it)}dt=-i\int_{\Gamma _{3}}\frac{1}{z}d\gamma=2\pi
n.
\]
Thus,
\begin{align*}
\sum_{k=1 }^Nz_k &=\frac{1}{
2\pi }\hspace{-0.05cm}\left[2\pi\lambda n-\pi \lambda N-2\pi i Nk  + 2\pi \lambda N+\pi i N^{2}\right]  \\
&=\frac{N\lambda }{2}+\lambda n+i\left(\frac{N^{2}}{2}+Nk\right),\quad k,n\in\mathbb{Z}.
\end{align*}
Finally, if $\Gamma$ contains at least one zero of $
\mathbf{B}_{(\lambda,N)}\varphi$, then there exists $z^\ast\in\mathbb{C}$ such that $z^\ast+\Gamma$ does not contain any zero as every nonzero analytic function can only have   finitely many zeros on any  compact set. The previous arguments can then be repeated for the curve $z^\ast+\Gamma$ yielding the same number of zeros which are constrained by the same condition. 
\pbox

\noindent The next result is a full characterization of frames obtained from the STFT on $S_N$ with Gaussian windows via sampling points in $\mathbb{T}_{N}^{2}$.

\begin{theorem}
\label{cor:frame-torus} Let $\lambda >0$, and   $z_{1},...,z_{K}\in \mathbb{T}_{N}^{2}$ be a collection
of $K$ distinct points. The following are equivalent: 
\begin{enumerate}\item[1).] The   points $z_{1},...,z_{K}\in \mathbb{T}_{N}^{2}$ give rise to
a frame for $S_{N}$, i.e., there are constants $A,B>0$ such that, for every $
\varphi \in S_{N}$,  
\begin{equation*}
A\Vert \varphi \Vert _{S_{N}}^{2}\leq \sum_{k=1}^{K}|\mathbf{V}
_{h_{0}^{\lambda }}\varphi (z_{k})|^{2}\leq B\Vert \varphi \Vert _{S_{N}}^{2}
\text{.}
\end{equation*}
\item[2).]
One of the following two conditions is satisfied:
\begin{enumerate}
    \item[(i)] $K\geq N+1$,
    \item[(ii)] $K=N$ and $\frac{1}{N}\sum_{k=1}^Nz_k\neq (1/2+n/N,N/2+m),$ for every $n,m\in\mathbb{Z}$.
\end{enumerate}
\end{enumerate}
\end{theorem}
\proof If the family is a frame for $S_N$, then it has to include at least $N$ vectors, i.e., $K\geq N$.

Let $K\geq N+1$. It follows from Proposition~\ref{thm:N-zeros} and the relation between $B_{(\lambda,N)}$ and $\mathbf{V}
_{h_{0}^{\lambda }}$ that $\mathbf{V}_{h_0^\lambda}\varphi $ has at most $N$ distinct zeros. Consequently, $\sum_{k=1}^{K}|\mathbf{V}
_{h_{0}^{\lambda }}\varphi (z_{k})|^{2}$ is always positive. Moreover, this expression
depends continuously on the basis coefficients $a_{n}$ of $\varphi$ by \eqref{eq:V_gvarphi}. Therefore, as the
unit sphere in $\mathbb{C}^{N}$ is compact, it follows that 
\begin{equation*}
A=\inf_{\Vert \varphi \Vert _{S_{N}}^{2}=1}\sum_{k=1}^{K}|\mathbf{V}
_{h_{0}^{\lambda }}\varphi (z_{k})|^{2}>0.
\end{equation*}
Now, let $K=N$, and suppose we are given a collection of $N$ distinct points $z_{1},...,z_{N}\in 
\mathbb{T}_{N}^{2}$. Moreover, we set $z_{\lambda,k}=\lambda x_k-i\xi_k\in\mathbb{T}^2_{\lambda,N}$, $k=1,...,N$. Let $z_{0}=\lambda/2+iN/2$ be the single zero of $
 \vartheta \left( iz/N,i\lambda/N \right) $ in $[0,\lambda]\times i[0,N]$ (see \cite[20.2(iv)]{NIST10}) and define the
function $F(z)$ as
\begin{equation*}
F(z)=\prod\limits_{k=1}^{N} e^{2\pi \text{Re}(z_{\lambda,k}-z_0) z/\lambda N} \vartheta \big(
i(z-z_{\lambda,k}+z_{0})/N,i\lambda/N \big) \text{.}
\end{equation*}
Using \eqref{eq:jacobi-quasi} one can directly show that $F$ satisfies the periodicity conditions
\begin{align*}
F(z+\lambda  n)&=
\prod\limits_{k=1}^{N} e^{2\pi \text{Re}(z_{\lambda,k}-z_0) (z+\lambda n)/\lambda N} \vartheta \big(
i(z-z_{\lambda,k}+z_{0})/N+i\lambda n/N,i\lambda/N \big)
\\
&=e^{\pi \lambda n^2}F(z)\prod\limits_{k=1}^{N} e^{2\pi \text{Re}(z_{\lambda,k}-z_0)   n/  N} e^{2\pi n(z-z_{\lambda,k}+z_0)/N} 
\\
&=e^{\pi \lambda n^2}e^{2\pi nz}e^{-2\pi in\text{Im}\big(Nz_0- \sum_{k=1}^Nz_{\lambda,k}\big)/ N}F(z),
\end{align*}
and
\begin{align*}
\quad F(z+iNm)&=\prod\limits_{k=1}^{N} e^{2\pi \text{Re}(z_{\lambda,k}-z_0) (z+N m)/\lambda N} \vartheta \big(
i(z-z_{\lambda,k}+z_{0})/N- m,i\lambda/N \big)
\\
&=e^{-2\pi im\text{Re}\big(Nz_0- \sum_{k=1}^Nz_{\lambda,k}\big)/ \lambda}F(z).
\end{align*}
Therefore, if  $\sum_{k=1}^Nz_{\lambda,k} = Nz_0+\lambda m+iNn=\lambda N/2+\lambda m+i(N^2/2+Nn),$ for some $n,m\in\mathbb{Z},$  it follows that $F$ satisfies 
\eqref{eq:barg-time-shift} and \eqref{eq:barg-freq-shift} and is
analytic in $\mathbb{T}_{\lambda,N}^{2}$. Thus, by Lemma~\ref
{dimension} there exists $\varphi \in S_{N}$
such that $\mathbf{V}_{h_{0}^{\lambda }}\varphi \left(\frac{x}{\lambda},-\xi\right)=F(z)e^{-\pi x^2/\lambda}$ and, by construction, $\mathbf{V}
_{h_{0}^{\lambda }}\varphi (z)=0$, for $z=z_{\lambda,1},...,z_{\lambda,K}$. Consequently,  
\begin{equation*}
\sum_{k=1}^{K}|\mathbf{V}
_{h_{0}^{\lambda }}\varphi (z_{\lambda,k})|^{2}=0,
\end{equation*}
and the lower frame bound is violated.

If $\sum_{k=1}^Nz_{\lambda,k}\neq Nz_0+\lambda m+Nn$, for every $n,m\in\mathbb{Z},$
it follows  by Proposition~\ref{thm:N-zeros} that $\{z_{\lambda,1},...,z_{\lambda,N}\}$ is a uniqueness set for the space of entire functions satisfying the periodicity conditions \eqref{eq:barg-time-shift} and \eqref{eq:barg-freq-shift}. It thus follows again by compactness of the unit ball in $S_N$, that the points $z_1,...,z_N$ generate a frame for $S_N$ 
\pbox

\noindent Theorem~\ref{cor:frame-torus} implies that if one chooses $N$ points in $\mathbb{T}^2_N$ uniformly at random, then one obtains a frame with probability 1. While the frame set for the Gaussian window STFT on $\mathrm{L}^2(\mathbb{R})$ is known from \cite{sewa92} and \cite{ly92}, it seems to be a folklore result, backed up by substantial evidence through numerical computations, that the finite Gabor transform with sampled, periodized Gaussian yields a frame when sampled on any lattice within $\mathbb{Z}_N\times \mathbb{Z}_N$ with cardinality larger than $N$, but we found no proof in the literature. On the other hand, a result by S\o{}ndergaard~\cite{sondergaard2007finite}, adapted from a remarkable observation by Benedetto et al~\cite{Benedetto1998} for the STFT on $\ell^2(\mathbb{Z})$, demonstrates that half-point shifts of sampled, periodized Gaussians yield frames on lattices of cardinality $N$. Results closer in spirit of our Theorem below, for totally positive functions and regular subsampling of $N\times N$, have been obtained in \cite[Theorem 9]%
{bannert2013discretized} using quite different methods.
The result below, a direct consequence of Theorem \ref{cor:frame-torus}, demonstrates that, in fact, any $K>N$ distinct points from the set $\{(j,k):\ j,k\in 0,...,N\}$ yield a finite Gabor frame for $\mathbb{C}^N$ with the sampled, periodized Gaussian window $\mathsf{h_N^\lambda}:=\mathbf{P}_{N}h_{0}^\lambda$. For $K=N$, additional arithmetic conditions ensure the frame property. This is the result for  Gabor frames in finite dimensions
stated in the introduction as one of the main achievements of our theory. We state it again for convenience.

\newcounter{orgthms}
\setcounter{orgthms}{\value{theorem}}
\setcounter{theorem}{2}
\begin{theorem} 
Let $\lambda>0$, and $\{(j_{k},l_{k})\}_{k=1,...,K}$ a collection of distinct pairs of integers  $j_{k},l_{k}\in 0,....,N-1$. The following are equivalent:
\begin{enumerate}
    \item[1).] 
The set $
\{(j_{k},l_{k})\}_{k=1,...,K}$    gives
rise to a finite Gabor frame with window $\mathsf{h_{N}^\lambda}$, i.e., there are constants $A,B>0$ such that, for
every $\mathbf{f}\in \mathbb{C}^N$,   
\begin{equation*}
A\left\|  \mathsf{f} \right\|_{\mathbb{C}^{N}}^{2}\leq
\sum_{k=1}^{K}\left\vert \V_{\mathsf{h_{N}^\lambda}}\mathsf{f}[j_k,l_k]\right\vert ^{2}\leq B \left\|\mathsf{f}\right\| _{\mathbb{C}^{N}}^{2}\text{.}
\end{equation*}
\item[2).] One of the three following conditions is satisfied:
\begin{enumerate}
    \item[(i)] $N^{2}\geq K\geq N+1,$
    \item[(ii)] $K=N$ is odd, 
\item[(iii)] $K=N$ is even and $\sum_{k=1}^N(j_k,l_k)\notin N\mathbb{N}^2. $
\end{enumerate}
\end{enumerate}
\end{theorem}
\setcounter{theorem}{\value{orgthms}}
\proof
Select a collection of $K$ distinct points $z_{1},...,z_{K}\in \mathbb{T}_{N}
$ of \ the form $\left(\frac{j_{k}}{N},l_{k}\right)\in \mathbb{T}_{N}$, $j_{k},l_{k}\in
0,....,N-1$ and
rewrite Theorem~\ref{cor:frame-torus} as 
\begin{equation*}
A\Vert \mathbf{\Sigma }_{N}f\Vert _{S_{N}}^{2}\leq \sum_{k=1}^{K}\left\vert
\mathbf{V}_{h_{0}^\lambda}\mathbf{\Sigma} _{N}f\left( \frac{j_{k}}{N},l_{k}\right) \right\vert
^{2}\leq B\Vert \mathbf{\Sigma }_{N}f\Vert _{S_{N}}^{2}\text{.}
\end{equation*}
Lemma~\ref{lem:double-period-2}  ensures that for every $\mathsf{f}\in\mathbb{C}^N$ there exists  $ {f}\in S_{0}(\mathbb{R})$ such that $\mathsf{P_N} {f}=\mathsf{f}$.
By \eqref{eq:a_n}
\begin{equation*}
\Vert \mathbf{\Sigma }_{N} {f}\Vert _{S_{N}}=\left\Vert \frac{1}{N}\mathsf{P_{N}} {f}\right\Vert _{\mathbb{C}^{N}}=\frac{1}{N}\|\mathsf{f}\|_{\mathbb{C}^N}.
\end{equation*}
Finally, by Theorem~\ref{thm:finite-gabor}, 
\begin{equation*}
\mathbf{V}_{h_{0}^\lambda}\mathbf{\Sigma} _{N}f\left( \frac{j_{k}}{N},l_{k}\right) =N^{-1}\mathbf{V}_{\mathsf{h_{N}^\lambda}}\mathsf{f}[j_k,l_k]. 
\end{equation*}
By Theorem~\ref{cor:frame-torus}, it follows that the frame conditions are satisfied if and only if either $K\geq N+1$, or $K=N$ and $\sum_{k=1}^N (j_k,l_k)\neq (N^2/2+Nn,N^2/2+Nm)$, for all $n,m\in\mathbb{Z}$. 

If $N=K$ is odd, then $N^2/2+nN=Nk+1/2$, for some $k\in\mathbb{Z}$. Therefore, the condition (ii) is automatically satisfied as $\sum_{k=1}^N(j_k,l_k)\in \mathbb{N}^2$.

If $N=K$ is even, then for every $k\in\mathbb{Z}$ there exist $n\in\mathbb{Z}$ such that $N^2/2+nN=Nk$. Therefore, we get a frame if and only if $\sum_{k=1}^N(j_k,l_k)\notin  N\mathbb{N}^2$.
\pbox

\noindent This result can be reformulated in the following sense: the discrete Gabor system, generated by the sampled and periodized Gaussian, is in \emph{general linear position} if $N$ is odd (and \emph{almost} in general linear position if $N$ is even). If any $N$ points from the set $\{(j,k):\ j,k\in 0,...,N\}$ yield a finite Gabor frame for $\mathbb{C}^N$ with window $\mathsf{g}$, then $\mathsf{g}$ is said to be in general linear position, see for example \cite{Goetz}. In \cite{Malikiosis} it is shown that almost every vector in $\mathbb{C}^N$ is in general linear position for every $N\in\mathbb{N}$. However, explicitly known examples of vectors in general linear position  are not localized and thus of no  use for practical purposes of finite Gabor analysis. Our result on the other hand allows to use localized windows (for appropriate choices of $\lambda$) for the prize that one needs to use one additional sampling point  if the number of points $N\in\mathbb{N}$ is even and enjoys a particular arithmetic structure.
{The property that a vector $\mathsf{g}\in\mathbb{C}^N$ is in general linear position is  equivalent to the 
Gabor matrix $\mathcal{G}(\mathsf{g})\in\mathbb{C}^{N\times N^2}$, $ \mathcal{G}(\mathsf{g} )_{n,(j,k)}=\mathsf{g} [n-j]e^{2\pi i kn/N}$, being \emph{full spark}.  We have thus shown that $\mathcal{G}(\mathsf{h_N^\lambda})$ is \emph{full spark if $N$ is odd}.}
The concept of the spark of a matrix has been considered by Donoho and Elad \cite{donoho2003optimally} as a measure to  decide uniqueness of the solution in sparse reconstruction problems. Later, Alexeev, Cahill, and Mixon introduced \emph{full spark frames}  \cite{alexeev2012full}. They proved that testing whether a given matrix is
full spark, is NP-hard under randomized polynomial-time reductions, concluding that \emph{deterministic full spark constructions} (such as the one we offer in this paper) \emph{are particularly significant}, because they
guarantee a property which is otherwise difficult to check. Full spark matrices/frames are of great importance in various applications, for example, in
compressed sensing \cite{fora13},   operator identification \cite{kalepfpo} and (Gabor) phase retrieval \cite{Palina}. 

Theorem \ref{cor:frame-torus} and Theorem \ref{cor:char-frame-set} can also be seen as 'Nyquist-type' necessary and sufficient result. The analogue necessary and sufficient result for infinite-dimensional Gabor frames with Gaussian window was proved by Seip and Wallst{\'e}n \cite{sewa92,se92-1} and independently by Lyubarskii \cite{ly92}. While the infinite dimensional results are stated in terms of Beurling density, in the toric setting we have a simple criteria: with less than $N$ samples the system is never a frame, but by increasing the number of samples above $N$, one is assured to have a frame with higher redundancy. Such a property comes in handy for situations where one is given a signal representation sampled on a grid with more than $N$ points (allowing for perfect reconstruction by our result) and wishes to increase the grid resolution for some numerical purpose. Then sampling again at a higher density is possible, and it still leads to perfect reconstruction. This can be done until all the possible $N^{2}$ points of the grid are used. 

By resorting to the extension to the torus and to Theorem \ref{cor:frame-torus} one can further use off-grid points and increase the resolution to arbitrary levels. 
Together with the first property of Proposition~\ref{thm:N-zeros}, where it is shown that the toric Bargmann transform of the signal has exactly $N$ zeros, this provides an ideal framework for the detection of spectrogram zeros, since knowing the number of zeros a priori is a clear advantage in such problems, in contrast with the infinite dimensional case, where this precise information is not available \cite{flandrin2015}. Besides, the constraint (\ref{eq:sum-zeros}) may be used to verify the precision of the detected zeros. It is also worth noticing that, while the
Hadamard-Weierstrass factorization of an entire function as an infinite
product cannot be used for reconstruction purposes, the situation is more
favorable for the toric Bargmann transform $\mathbf{B}_{(\lambda ,N)}\varphi 
$, since it can be written as the following finite product of $N$ theta
functions, each of them vanishing on one of their $N$ zeros $\{z_{\lambda
,k}\}_{k=1}^{N}$:
\begin{equation*}
\mathbf{B}_{(\lambda ,N)}\varphi (z)=C\prod\limits_{k=1}^{N}e^{2\pi \text{Re}
(z_{\lambda ,k}-z_{0})z/\lambda N}\vartheta \left( i(z-z_{\lambda
,k}+z_{0})/N,i\lambda /N\right) \text{,}
\end{equation*}
for some normalization constant $C$, and where $z_{0}=\lambda /2+iN/2$ is the single
zero of $\vartheta \left( iz/N,i\lambda /N\right) $ in $[0,\lambda ]\times
i[0,N]$ (this can be proved with the methods of \cite[Proposition 4.1]
{AFTALION2006661}). This allows to effectively parametrize  $
\mathbf{B}_{(\lambda ,N)}\varphi $ by its set of $N$ zeros, in a one to
one correspondence known as \emph{the stellar representation }\cite
{nonnenmacher1998chaotic}. 

Finally, we point out the model of \cite
{bardenet2020zeros}, where white noise, defined as a random linear
combination of Hermite functions with i.i.d. coefficients, is mapped by the
Bargmann transform into a Gaussian Entire Function, whose known zero
statistics are leveraged for the problem of reconstruction of a signal
embedded in white noise from the spectrogram zeros of the mixture, a method
put forward in \cite{flandrin2015}. In our setting, white noise $\mathcal{W}$
can be defined in an analogue way in $S_{N}$, using the basis \ref{def:deltatrain}, and circumventing the technicalities of the infinite dimensional model, by simply considering 
 $
\mathcal{W}=\sum_{n=0}^{N-1}a_{n}\epsilon _{n}, 
$ 
where $a_{n}$ are i.i.d. Gaussian random variables. Applying the Bargmann-type transform with $\lambda=1$

yields the following Gaussian Analytic Function on the flat tori
\begin{equation*}
\mathbf{B}_{(1,N)}\mathcal{W}(z)=\sum_{n=0}^{N-1}a_{n}e^{-\pi \left( \frac{n}{N}\right) ^{2}}e^{2\pi iz
\frac{n}{N}}\vartheta \left( i\left( z-\frac{n}{N}\right) ,i\right).
\end{equation*}
Some statistics of the distribution of the zeros of this Gaussian function (with
different normalizations) have been computed in \cite[Section 5.1.1]
{nonnenmacher1998chaotic}.

\bibliographystyle{plain}
\bibliography{project}

\begin{thebibliography}{10}

\bibitem{AFTALION2006661}
A.~Aftalion, X.~Blanc, and F.~Nier.
\newblock Lowest {L}andau level functional and {B}argmann spaces for
  {B}ose–{E}instein condensates.
\newblock {\em J. Funct. Anal.}, 241(2):661--702, 2006.

\bibitem{ahl}
L.~Ahlfors.
\newblock {\em Complex {A}nalysis: {A}n {I}ntroduction to the {T}heory of
  {A}nalytic {F}unctions of one {C}omplex {V}ariable}.
\newblock McGraw-Hill, 1979.

\bibitem{alexeev2012full}
B.~Alexeev, J.~Cahill, and D.~G. Mixon.
\newblock Full spark frames.
\newblock {\em J. Fourier Anal. Appl.}, 18:1167--1194, 2012.

\bibitem{Ausl96}
L.~Auslander and Y.~Meyer.
\newblock A generalized {P}oisson summation formula.
\newblock {\em Appl. Comp. Harmon. Anal.}, 3(4):372--376, 1996.

\bibitem{bannert2013discretized}
S.~Bannert, K.~Gr{\"o}chenig, and J.~St{\"o}ckler.
\newblock Discretized {G}abor frames of totally positive functions.
\newblock {\em IEEE Trans. Inform. Theory}, 60(1):159--169, 2013.

\bibitem{bardenet2020zeros}
R.~Bardenet, J.~Flamant, and P.~Chainais.
\newblock On the zeros of the spectrogram of white noise.
\newblock {\em Appl. Comp. Harmon. Anal.}, 48(2):682--705, 2020.

\bibitem{Benedetto1998}
J.~J. Benedetto, C.~Heil, and D.~F. Walnut.
\newblock Gabor systems and the {B}alian-{L}ow theorem.
\newblock In H.~G. Feichtinger and T.~Strohmer, editors, {\em Gabor Analysis
  and Algorithms: Theory and Applications}, pages 85--122. Birkh{\"a}user
  Boston, 1998.

\bibitem{coro}
M.~Combescure and D.~Robert.
\newblock {\em Coherent States and Applications in Mathematical Physics}.
\newblock Springer, 2012.

\bibitem{Cool67}
J.~Cooley, P.~Lewis, and P.~Welch.
\newblock Application of the fast {F}ourier transform to computation of
  {F}ourier integrals, {F}ourier series, and convolution integrals.
\newblock {\em IEEE Trans. Audio and Electroacoustics}, 15(2):79--84, 1967.

\bibitem{courbot2023sparse}
J.-B. Courbot, A.~Moukadem, B.~Colicchio, and A.~Dieterlen.
\newblock Sparse off-the-grid computation of the zeros of {STFT}.
\newblock {\em IEEE Signal Process. Lett.}, 2023.

\bibitem{donoho2003optimally}
D.~L. Donoho and M.~Elad.
\newblock Optimally sparse representation in general (nonorthogonal)
  dictionaries via $\ell^1$ minimization.
\newblock {\em Proc. Nat. Acad. Sci.}, 100(5):2197--2202, 2003.

\bibitem{escudero2022efficient}
L.~A. Escudero, N.~Feldheim, G.~Koliander, and J.~L. Romero.
\newblock Efficient computation of the zeros of the {B}argmann transform under
  additive white noise.
\newblock {\em Found. Comp. Math.}, pages 1--34, 2022.

\bibitem{FeichingerSegal}
H.~G. Feichtinger.
\newblock On a new {S}egal algebra.
\newblock {\em Monatsh. Math.}, 92:269--289, 1981.

\bibitem{fe06}
H.~G. Feichtinger.
\newblock Modulation spaces: Looking back and ahead.
\newblock {\em Samp. Theory Signal Image Process.}, 5:109--140, 2006.

\bibitem{FeGr1}
{H}.~{G}. {F}eichtinger and {K}. {G}r{\"o}chenig.
\newblock {B}anach spaces related to integrable group representations and their
  atomic decompositions, {I}.
\newblock {\em J. Funct. Anal.}, 86(2):307--340, 1989.

\bibitem{flandrin2015}
P.~Flandrin.
\newblock Time–frequency filtering based on spectrogram zeros.
\newblock {\em IEEE Signal Process. Lett.}, 22(11):2137--2141, 2015.

\bibitem{fora13}
S.~Foucart and H.~Rauhut.
\newblock {\em A {M}athematical {I}ntroduction to {C}ompressive {S}ensing}.
\newblock Birkh{\"a}user, 2013.

\bibitem{gardner2006sparse}
T.~J. Gardner and M.~O. Magnasco.
\newblock Sparse time-frequency representations.
\newblock {\em Proc, Nat. Acad. Sci.}, 103(16):6094--6099, 2006.

\bibitem{Gazeau}
J.-P. Gazeau.
\newblock {\em Coherent {S}tates in {Q}uantum {P}hysics}.
\newblock Wiley, 2009.

\bibitem{Charly}
K.~Gr{\"o}chenig.
\newblock {\em {F}oundations of {T}ime-{F}requency {A}nalysis}.
\newblock {A}ppl. {N}umer. {H}armon. {A}nal. {B}irkh{\"a}user {B}oston, 2001.

\bibitem{jakobsen}
M.~S. Jakobsen.
\newblock On a (no longer) new {S}egal algebra.
\newblock {\em J. Fourier Anal. Appl.}, 24(6):1579--1660, 2018.

\bibitem{jalu19}
M.~S. Jakobsen and F.~Luef.
\newblock Sampling and periodization of generators of {H}eisenberg modules.
\newblock {\em Int. J. Math.}, 30(10):1950051, 2019.

\bibitem{Jans97}
A.~J. E.~M. Janssen.
\newblock From continous to discrete {W}eyl-{H}eisenberg frames through
  sampling.
\newblock {\em J. Fourier Anal. Appl.}, 3(5):583--596, 1997.

\bibitem{Kaib05}
N.~Kaiblinger.
\newblock Approximation of the {F}ourier transform and the dual {G}abor window.
\newblock {\em J. Fourier Anal. Appl.}, 11:25--42, 02 2005.

\bibitem{kalepfpo}
A.~Kaplan, D.~Lee, G.~E. Pfander, and V.~Pohl.
\newblock Sparse deterministic and stochastic channels: Identification of
  spreading functions and covariances.
\newblock In G.~Kutyniok, H.~Rauhut, and R.~J. Kunsch, editors, {\em Compressed
  {S}ensing in {I}nformation {P}rocessing}, chapter~4. Birkh{\"a}user Boston,
  2022.

\bibitem{katz}
Y.~Katznelson.
\newblock {\em An {I}ntroduction to {H}armonic {A}nalysis}.
\newblock Cambridge University Press, 3rd edition, 2004.

\bibitem{LV}
P.~Leboeuf and A.~Voros.
\newblock Chaos-revealing multiplicative representation of quantum eigenstates.
\newblock {\em J. Phys. A: Mathematical and General}, 23(10):1765, 1990.

\bibitem{ly92}
{Y}.~{I}. {L}yubarskii.
\newblock {F}rames in the {B}argmann space of entire functions.
\newblock In {\em {E}ntire and {S}ubharmonic {F}unctions}, volume~11 of {\em
  {A}dv. {S}ov. {M}ath.}, pages 167--180. {A}merican {M}athematical {S}ociety
  ({A}{M}{S}), 1992.

\bibitem{Malikiosis}
R.~D. Malikiosis.
\newblock A note on {G}abor frames in finite dimensions.
\newblock {\em Appl. Comp. Harmon. Anal.}, 38(2):318--330, 2015.

\bibitem{nonnenmacher1998chaotic}
S.~Nonnenmacher and A.~Voros.
\newblock Chaotic eigenfunctions in phase space.
\newblock {\em J. Stat. Phys.}, 92:431--518, 1998.

\bibitem{NIST10}
F.~W. Olver, D.~W. Lozier, R.~F. Boisvert, and C.~W. Clark.
\newblock {\em {NIST} {H}andbook of {M}athematical {F}unctions}.
\newblock Cambridge University Press, 1st edition, 2010.

\bibitem{Goetz}
G.~E. Pfander.
\newblock Gabor frames in finite dimensions.
\newblock In P.~G. Casazza and G.~Kutyniok, editors, {\em Finite Frames: Theory
  and Applications}, pages 193--239. Birkh{\"a}user Boston, Boston, 2013.

\bibitem{Palina}
P.~Salanevich and G.~E. Pfander.
\newblock Geometric properties of {G}abor frames with a random window.
\newblock In {\em 2017 International Conference on Sampling Theory and
  Applications (SampTA)}, pages 183--187, 2017.

\bibitem{se92-1}
{K}. {S}eip.
\newblock {D}ensity theorems for sampling and interpolation in the
  {B}argmann-{F}ock space. {I}.
\newblock {\em {J}. {R}eine {A}ngew. {M}ath.}, 429:91--106, 1992.

\bibitem{sewa92}
{K}. {S}eip and {R}. {W}allst{\'e}n.
\newblock {D}ensity theorems for sampling and interpolation in the
  {B}argmann-{F}ock space. {I}{I}.
\newblock {\em {J}. {R}eine {A}ngew. {M}ath.}, 429:107--113, 1992.

\bibitem{sondergaard2007finite}
P.~S{\o}ndergaard.
\newblock {\em Finite discrete {G}abor analysis}.
\newblock PhD thesis, Institut for Matematik, DTU, 2007.

\bibitem{so05}
P.~Søndergaard.
\newblock Gabor frames by sampling and periodization.
\newblock {\em Adv. Comput. Math.}, 27:355--373, 10 2007.

\end{thebibliography}

\end{document}